\documentclass[12pt]{article}

\usepackage[english]{babel}
\usepackage{latexsym, graphicx, natbib}

\usepackage{amsthm,amsmath,natbib}
\RequirePackage[dvips]{hyperref}

\DeclareSymbolFont{AMSb}{U}{msb}{m}{n}
\DeclareSymbolFontAlphabet{\mathbb}{AMSb}

\newcommand{\nc}{\newcommand}
\nc{\nt}{\newtheorem} \nt{defn}{Definition} \nt{lem}{Lemma}
\nt{pr}{Proposition}
\nt{theorem}{Theorem} \nt{cor}{Corollary} \nt{ex}{Example}
\nt{ass}{Assumption} \nt{step}{Step} \nt{case}{Case}
\nt{subcase}{Subcase} \nt{note}{Remark} \nc{\bd}{\begin{defn}}
\nc{\ed}{\end{defn}} \nc{\blem}{\begin{lem}} \nc{\elem}{\end{lem}}
\nc{\bpr}{\begin{pr}} \nc{\epr}{\end{pr}} \nc{\bth}{\begin{theorem}}
\nc{\eth}{\end{theorem}} \nc{\bcor}{\begin{cor}}
\nc{\ecor}{\end{cor}} \nc{\bex}{\begin{ex}}  \nc{\eex}{\end{ex}}
\nc{\bass}{\begin{ass}}  \nc{\eass}{\end{ass}}
\nc{\bstep}{\begin{step}}  \nc{\estep}{\end{step}}
\nc{\bcase}{\begin{case}}  \nc{\ecase}{\end{case}}
\nc{\bsubcase}{\begin{subcase}}  \nc{\esubcase}{\end{subcase}}
\nc{\bnote}{\begin{note}}  \nc{\enote}{\end{note}} \nc{\prf}{{\bf
Proof} }
\nc{\eop}{\hfill $\Box$ \\ \\}%\Binom{n}{k}
\nc{\argmin}{\mathrm{argmin}} \nc{\argmax}{\mathrm{argmax}}
\nc{\sgn}{\mathrm{sgn}} \nc{\supp}{\mathrm{supp}}
\nc{\Var}{\mathrm{Var}} \nc{\Cov}{\mathrm{Cov}}
\nc{\bak}{\!\!\!\!\!} \nc{\IBD}{\mathrm{IBD}}

\begin{document}

\title{Limit properties of the monotone rearrangement for density and regression function estimation}
\author{Dragi Anevski and Anne-Laure Foug\`eres\\ G\"oteborg University and Universit\'e Paris X - Nanterre}
\maketitle

\begin{abstract}
~~~~The monotone rearrrangement algorithm was introduced by Hardy,
Littlewood and P\'olya as a sorting device for functions.  Assuming
that $x$ is a monotone function and that an estimate $x_n$ of $x$ is
given, consider the monotone rearrangement $\hat{x}_n$ of $x_n$.
This new estimator is shown to be uniformly consistent. Under
suitable assumptions, pointwise limit distribution results for
$\hat{x}_n$ are obtained. The framework is general and allows for
weakly dependent and long range dependent stationary data.
Applications in monotone density and regression function estimation
are detailed.
\end{abstract}


\section{Introduction}
Assume that $(t_i,x(t_i))_{i=1}^n$, for some points $t_i\in [0,1]$
 (e.g. ($t_i=i/n$)), are pairs of data points.  The (decreasing)
sorting of the points $x(t_i)$ is then an elementary operation and
produces the new sorted sequence of pairs $(t_i,y(t_i))$ where
$y=sort(x)$ is the sorted vector. Let $\#$ denote the counting
measure of a set. Then we can define the sorting $y$ of $x$ by
\begin{eqnarray*}
       z(s)&=&\#\{t_i:x(t_i)\geq s\}\\
       y(t)&=&z^{-1}(t),
\end{eqnarray*}
where $z^{-1}$ denotes the inverse of a function (if the points
$x(t_i)$ are not unique it denotes the generalized inverse).  \par
The "sorting" of a function $\{x(t),t\in [0,1]\}$ can then
analogously be defined by the monotone rearrangement (cf. Hardy et
al. \cite{hardy:littlewood:polya:1952}),
\begin{eqnarray*}
     z(s)&=&\lambda\{t\in [0,1]:x(t)\geq s\},\\
     y(t)&=&z^{-1}(t),
\end{eqnarray*}
where the counting measure $\#$ has been replaced by the Lebesgue
measure $\lambda$, and $z^{-1}$ denotes the generalized inverse.
\par The monotone rearrangement algorithm of a set or a function
has mainly been used as a device in analysis, see e.g. Lieb and Loss
\cite[Chapter 3]{lieb:loss:1996} or in optimal transportation (see
Villani \cite[Chapter 3]{villani:2003}). Foug\`eres
\cite{fougeres:1997} was the first to use the algorithm in a
statistical context, for density estimation under order
restrictions. Meanwhile, Polonik \cite{polonik:1995, polonik:1998}
also developed tools of a similar kind for density estimation for
multivariate data. More recently, several authors revisited the
monotone rearrangement procedure in the estimation context under
monotonicity; see Dette et al. \cite{dette:neumeyer:pilz:2006}, and
Chernozhukov et al. \cite{chernozhukov:fernandezval:galichon:2007a}.
\par We introduce the following two-step approach for estimating
a monotone function. Assume that $x$ is a monotone function on an interval
$I \subset {\mathbb R}$. Assume also that we already have an
estimate $x_n$ of $x$, but that this estimate is not necessarily
monotone. We then propose to use the monotone rearrangement
$\hat{x}_n$ of $x_n$ as an estimate of $x$.
\par Under the assumption that we have process limit distribution
results for (a localized version of) the stochastic part of $x_n$
and that the deterministic part of $x_n$ is asymptotically
differentiable at a fixed point $t_0$, with strictly negative
derivative, we obtain pointwise limit distribution results for
$\hat{x}_n(t_0)$. The framework is general and allows for weakly
dependent as well as long range dependent data. This is the topic
for Section 3.  \par Possible applications of the general results
are to monotone density and regression function estimation, which we
explore in more detail in Section 4. These are the problems of
estimating $f$ and $m$ respectively in
\begin{eqnarray*}
      (i)&&t_1,\ldots,t_n\mbox{ stationary observations with marginal}\\
      &&\mbox{decreasing density } f \mbox{ on }{\mathbb R^{+}},\\
       (ii)&&(t_i,y_i)\mbox{ observations from  } y_i=m(t_i)+\epsilon_i,\\
    &&t_i=i/n,i=1,\ldots,n, m\mbox{ decreasing on }[0,1],\\
    &&\{\epsilon_i\}\mbox{ stationary sequence with mean zero}.
\end{eqnarray*}
The standard approaches in these two problems have been isotonic
regression for the regression problem, first studied by Brunk
\cite{brunk:1955}, and (nonparametric) Maximum Likelihood estimation
(NPMLE) for the density estimation problem, first introduced by
Grenander \cite{grenander:1956}. A wide literature exists for
regression and density estimation under order restrictions. One can
refer e.g. to Mukerjee \cite{mukerjee:1988}, Ramsay
\cite{ramsay:1988}, Mammen \cite{mammen:1991}, Hall and Huang
\cite{hall:huang:2001}, Mammen et al.
\cite{mammen:marron:turlach:wand:2001}, Gijbels \cite{gijbels:2004},
Birke and Dette \cite{birke:dette:2007}, Dette and Pilz
\cite{dette:pilz:2006}, Dette et al. \cite{dette:neumeyer:pilz:2006}
for the regression context. Besides, see Eggermont and Lariccia
\cite{eggermont:lariccia:2000}, Foug\`eres \cite{fougeres:1997},
Hall and Kang \cite{hall:kang:2005}, Meyer and Woodroofe
\cite{meyer:woodroofe:2004}, Polonik \cite{polonik:1995}, Van der
Vaart and Van der Laan \cite{vandervaart:vanderlaan:2003}, among
others, for a focus on monotone (or unimodal) density estimation.
Anevski and H\"ossjer \cite{anevski:hossjer:2006} gave a general
approach unifying both contexts.

 Using kernel estimators as preliminary estimators of $f$ and $m$ on which the
monotone rearrangement is then applied, we are able to derive limit
distribution results for quite general dependence situations,
demanding essentially stationarity for the underlying random parts
$\{t_i\}$ and $\{\epsilon_i\}$  respectively. The results are
however stated in a form that allows for other estimators than the
kernel based as starting points, e.g. wavelet or splines estimators.
\par
The paper is organized as follows: In Section 2 we define the
monotone rearrangement algorithm and derive some simple properties
that will be used in the sequel. In particular our definition
differs slightly from Hardy, Littlewood and Polya's original
definition \cite{hardy:littlewood:polya:1952}; the difference is
motivated by the fact that we will use localization and restriction.
The most important properties for the algorithm that are derived are
the equivariance under addition of constants, the continuity of the
map and a certain localization property, cf. Lemma \ref{lem:proprea},
Theorem \ref{thm:contraction} and Theorem \ref{thm:trunc}
below. Furthermore we state conditions that allow for the extension
of the map to unbounded intervals.
\par
In Section 3 we define the generic estimator of the monotone function, and
state the consistency and limit distribution properties for the estimator.
The limit distribution is given in Theorem 4 and is of the general form
\begin{eqnarray*}
      d_n^{-1}\left[  \hat{x}_{n}(t_0)-x(t_0)\right]&\stackrel{\cal L}{\rightarrow}&
       T\left(A \cdot +\tilde{v}(\cdot ;t_0)\right)(0)+\Delta,
\end{eqnarray*}
where $T$ is the monotone rearrangement map, $\Delta=\lim_{n\to
\infty} d_n^{-1} [{\mathbb E}\{x_n(t_0+sd_n)\}-x(t_0)]$ is the
asymptotic local bias of the preliminary estimator and
$\tilde{v}(s;t_0)\stackrel{\cal L}{=}\lim_{n\to\infty}d_n^{-1}
\left[x_n(t_0+sd_n)-{\mathbb E}\{x_n(t_0+sd_n)\}\right]$ is the weak
local limit of the process part of the preliminary estimator; here
$d_n\downarrow 0$ is a deterministic sequence that is  determined by
the dependence structure of the data.
\par
In Section 4 we apply the obtained results in Section 3 to
regression function estimation and density estimation under order
restrictions, and derive the limit distributions for the estimators.
This gives rise to some new universal limit random variables, such
as e.g. in the regression context $T(s+B(s))(0)$ with $T$ the
monotone rearrangement map and $B$ standard two sided Brownian
motion for independent and weakly dependent data, or
$T(s+B_{1,\beta}(s))(0)$ with $B_{1,\beta}$ fractional Brownian
motion with self similarity parameter $\beta$, when data are long
range dependent. The rate of convergence $d_n$ is e.g. for the
regression problem the optimal $n^{-1/3}$ in the i.i.d. and weakly
dependent data context and of a non-polynomial rate in the long
range dependent context, similarly to previously obtained results in
isotonic regression for long range dependent data, cf. Anevski and
H\"ossjer \cite{anevski:hossjer:2006}.
\par
In the appendix we derive some useful but technical results on
maximal bounds on the rescaled process parts in the density and
regression estimation problems, i.e. for the local partial sum
process and empirical processes, for weakly dependent as well as
long range dependent data.

\section{The monotone rearrangement algorithm}
Consider an interval $I\subset{\mathbb R}$, and let ${\cal
B}(I)=\{f: \, f(I) \mbox{ bounded}\}$ and ${\cal D}(I)=\{f:\,
f\mbox{ decreasing on }I\}$. For each Borel set $A$ of ${\mathbb
R}$, denote
 $\lambda(A)$ the Lebesgue measure of $A$ on ${\mathbb R}$. In a
 first step, the monotone rearrangement will be defined for finite intervals
 $I$, and some extensions for infinite $I$ will be discussed in a
 second step.
\subsection{Definition and properties for finite intervals}
 \bd Let $I\subset{\mathbb R}$ be a finite interval, and assume $f\in
{\cal B}(I)$. Let $r_{f,I}$ be the right continuous map from $f(I)$ to
${\mathbb R}^+$, called "upper level set function" of $f$ and
defined for each $u\in f(I)$ by
\begin{eqnarray*}
     r_{f,I}(u)&:=&\lambda\{t\in I :\, f(t)>u\}=\lambda\{I \cap f^{-1}(u,\infty)\}.
\end{eqnarray*}
The monotone rearrangement map $T_I:{\cal B}(I)\ni f\mapsto
T_I(f)\in{\cal D}(I)$ is defined up to a translation as the (right
continuous) generalized inverse of the upper level set function
   \begin{eqnarray}
 T_I(f)(t)&:=&\inf\{u\in f(I):r_{f,I}(u)\leq t- \inf I\}, \label{eq:T_I_def}
   \end{eqnarray}
for $t \in I$. \ed

The following lemmas are listing some simple and useful properties
of the maps $u\mapsto r_{f,I}(u)$, $f \mapsto r_{f,I}$ and $f \mapsto
T_I(f)$ respectively.

 \blem Assume $I\subset{\mathbb R}$ is a finite interval, and $f\in
{\cal B}(I)$. Then
\begin{eqnarray*}
    (i)&& \mbox{If $f$ has no flat regions on $I$, i.e.  $\lambda\{I\cap f^{-1}(u)\}=0$
    }\\
    &&\mbox{for all $u \in f(I)$, then $ r_{f,I}$ is continuous},\\
    (ii)&& \mbox{If there is a $u_{0} \in f(I)$ such that $\lambda\{I\cap f^{-1}(u_0)\}=c>0$
    then $r_{f,I}$}\\
    &&\mbox{has a discontinuity at $u_0$ of height c},\\
    (iii)&& \mbox{If $f$ has a discontinuity at $t_0 \in I$ and $f$ is decreasing,
    then $r_{f,I}$}\\
    &&\mbox{admits a flat region with level $t_0$}.
\end{eqnarray*}
\elem \prf Assertions $(i)$ and $(ii)$ are both consequences of the
fact that
\begin{eqnarray*}
\lim_{u\to u_0} |r_{f,I}(u)-r_{f,I}(u_0)| &=& \lim_{u \to u_0} \lambda \{ t
\in I: \max(u,u_0) \geq f(t) > \min(u,u_0) \} \\
&=& \lambda \{ I\cap f^{-1}(u_0)\},
\end{eqnarray*}
 which is equal to 0 in $(i)$,
and to $c$ in $(ii)$. Finally, assertion $(iii)$ arises from writing
that  $r_{f,I}(u)= r_{f,I}(f(t_0^-))=t_0$ for each $u \in
(f(t_0^+),f(t_0^-))$.
 \eop

\blem Let $I\subset{\mathbb R}$ be a finite interval, and assume
$f\in {\cal B}(I)$. Then
\begin{eqnarray*}
    (i)&&\mbox{If $c$ is a constant then } r_{f+c,I}(u)=r_{f,I}(u-c), \mbox{for each } u \in f(I)+c.\\
    (ii)&&r_{cf,I}(u)=r_{f,I}(u/c)\mbox{ if }c>0, \mbox{for each } u \in c f(I).\\
    (iii)&&f\leq g \Rightarrow r_{f,I}\leq r_{g,I}.\\
    (iv)&&\mbox{Let }f_c(t)=f(tc). \mbox{ Then }c \, r_{f_c,I}=r_{f,I}.\\
    (v)&&\mbox{Let }f_c(t)=f(t+c). \mbox{ Then }r_{f_{c},I}=r_{f,I}.
\end{eqnarray*}
\elem \prf (i)-(iii) follow from the definition; indeed, for each $u
\in f(I)+c$, $r_{f+c,I}(u)=\lambda\{t\in I:f(t)+c>u\}=r_{f,I}(u-c)$, and
for each $u \in c f(I)$, $r_{cf,I}(u)=\lambda\{t\in
I:cf(t)>u\}=r_{f,I}(u/c)$ if $c>0$. As for (iii), $\{t\in
I:f(t)>u\}\subset\{t\in I:g(t)>u\}$, for each fixed $u$, if $f\leq
g$. Statement (iv) follows from $r_{f_c,I}(u)=\lambda\{t\in I/c
:f(ct)>u\}=\lambda\{s/c \in I/c : f(s)>u\}=r_{f,I}(u)/c$, for each $u
\in f(I)$. Statement (v) is a consequence of
$r_{f_c,I}(u)=\lambda\{t\in I-c:f(t+c)>u\}=\lambda\{s-c \in
I-c:f(s)>u\}=\lambda\{t\in I :f(t)>u\}$, for each $u \in f(I)$. \eop
 \blem\label{lem:proprea} Let
$I\subset{\mathbb R}$ be a finite interval and assume $f,g$ are
functions in ${\cal B}(I)$. The monotone rearrangement map $T_I$
satisfies the following:
\begin{eqnarray*}
   (i)&&T_I(f+c)=T_I(f)+c, \mbox{ if $c$ is a constant};\\
   (ii)&&T_I(cf)=cT_I(f), \mbox{ if $c>0$ is a constant};\\
   (iii)&&f\leq g\Rightarrow T_I(f)\leq T_I(g);\\
   (iv)&&\mbox{Let }f_c(t)=f(ct); \mbox{ then }T_{I/c}(f_c)(t)=T_{I}(f)(ct);\\
   (v)&&\mbox{Let }f_{c}(t)=f(t+c); \mbox{ then }T_{I-c}(f_c)(t)=T_{I}(f)(t+c).
\end{eqnarray*}
\elem \prf Let $I = [a,b]$; each assertion is a consequence of its
counterpart in Lemma 2. Let $t \in I$; statement $(i)$ follows from
$T_I(f+c)(t)= \inf \{u \in f(I)+c : r_{f,I}(u-c) \leq t-a\} =
T_I(f)(t)+c$, whereas $(ii)$ comes from $T_I(cf)(t)= \inf\{u \in c
f(I) : r_{f,I}(u/c) \leq t-a \} = c T_I(f)(t)$. To show $(iii)$, note
that $f\leq g \Rightarrow r_{f,I}\leq r_{g,I} \Rightarrow T_I(f)\leq
T_I(g)$. Assertion $(iv)$ follows from the fact that for each $t \in
I/c$,  $T_{I/c}(f_c)(t)= \inf\{u \in f(I): r_{f,I}(u) \leq ct-a\} =
T_I(f)(ct)$. Finally, statement $(v)$ follows since for each $t \in
I-c$, $ T_{I-c}(f_c)(t) = \inf \{ u \in f(I) : r_{f,I}(u) \leq t+c -a
\}= T_I(f)(t+c)$.
 \eop
 The previous result implies that the map $T_I$ is continuous, as
 stated in the following theorem.
\bth\label{thm:contraction} Let $||\cdot||$ be an arbitrary norm on
${\cal B}(I)$. Then the map $T_I$ is a contraction, i.e.
$||T_I(f)-T_I(g)||\leq ||f-g||$. In particular, $T_I$ is a continuous  map,
i.e. for all $f_n,f\in{\cal B}(I)$,
\begin{eqnarray*}
  ||f_n-f||\rightarrow 0 &\Rightarrow& ||T_I(f_n)-T_I(f)||\rightarrow
  0,
\end{eqnarray*}
as $n$ tends to infinity. \eth \prf Let $f,g$ be functions in ${\cal
B}(I)$. Clearly $g(u)-||f-g||\leq f(u) \leq g(u)+||f-g||$, which by Lemma 3
$(i)$ and $(iii)$ implies that $T_I(g)(u)-||f-g||\leq T_I(f)(u) \leq
T_I(g)(u)+||f-g||$, so that $|T_I(f)-T_I(g)|(u)\leq ||f-g||$, for
each $u$. Since the right
hand side is independent of $u$, the absolute value on the left hand
side can be replaced by the norm, which implies the statement of the
theorem.  \eop
\bnote One can also refer to Lieb and Loss \cite[Theorem
3.5]{lieb:loss:1996}  for a proof
 of the contraction property (the "non expansivity" property of the map $T_I$), for the $L^p$-norms.
 \eop
\enote
\subsection{Extension to infinite intervals}
It is not possible to define a monotone rearrangement on an infinite
interval $I$ or on ${\mathbb R}$ for {\it any} function $\varphi \in
{\cal B}({\mathbb R})$. This can however be done for positive
functions $f$ such that for each $u>0$, $r_f(u):= \lambda\{t \in
{\mathbb R} : f(t)>u\} < + \infty$, defining in this situation
$$
T(f)(t) := \inf\{u\in {\mathbb R}_+ :r_f(u)\leq t  \},
$$
for each positive $t$. Such a definition is precisely the definition
considered by Hardy et al. \cite[Chapter
10.12]{hardy:littlewood:polya:1952}, and is in particular valid for
densities $f \in {\cal B}({\mathbb R})$ (see also Lieb and Loss
\cite[Chapter 3]{lieb:loss:1996} and Foug\`eres
\cite{fougeres:1997}).

If it remains impossible to define $T(\varphi)$ for any function
$\varphi$ for which $r_{\varphi}(u)$ is possibly infinite for some
positive $u$, such a definition can be given {\it locally} around a
fixed point $x \in I_0$, where $I_0$ is a finite interval, as soon as the
function $\varphi$ satisfies the following property:
\vskip 3mm \hskip 1.5cm There exists a constant $M<\infty$ and a
finite interval $I_1$ including

\hskip 1.5cm $I_0$  such that
\begin{eqnarray}
\inf_{t\in (\inf I_1,\sup I_0)} \varphi(t)  \geq -M &{\rm ~~and~~}&
\sup_{t\in (\sup I_1,\infty)} \varphi(t) \leq -M,\label{eq:trunc-hyp1}\\
\inf_{t\in (-\infty,\inf I_1)} \varphi(t)  \geq +M &{\rm ~~and~~}&
\sup_{t\in (\inf I_0,\sup I_1)} \varphi(t) \leq +M.\label{eq:trunc-hyp2}
   \end{eqnarray}

\bth\label{thm:trunc} Let $I_0$ be a finite and fixed interval, and
let $\varphi\in  {\cal B}({\mathbb R})$ such that
$(\ref{eq:trunc-hyp1})$ and $(\ref{eq:trunc-hyp2})$ are satisfied.
Then for any finite interval $J$ containing $I_1$, one has
$T_{J}(\varphi)\equiv T_{I_1}(\varphi)$ on $I_0$.  \eth
\prf Define $y_1:= \inf \{y \in I_1 : \forall x \in [\inf J, y[ \,
\; \varphi(x) > \varphi(y)\}$ and $z_0:= \inf\{x \in J : \varphi(x)
\in \varphi(I_0)\}$. If follows from those definitions, from the
left part of (\ref{eq:trunc-hyp2}) and from the continuity of
$\varphi$ that $y_1\in I_1$, $z_0 \in I_1$ and $y_1 <z_0 \leq \inf
I_0$. As a consequence, one has:
\begin{eqnarray*}
r_{\varphi,J}\{\varphi(y_1)\}& = & \lambda\{t\in J:\varphi(t)>
\varphi(y_1)\}\\
%& =& \lambda\{t\in J\cap(-\infty,y_1) :\varphi(t)> \varphi(y_1)\} +
%\lambda\{t\in J\cap(y_1,\infty) :\varphi(t)> \varphi(y_1)\}\\
&=& y_1 - \inf J + \lambda\{t\in I_1\cap(y_1,\infty) :\varphi(t)>
\varphi(y_1)\},
\end{eqnarray*}
where the second equality comes from splitting $J$ into
$J\cap(-\infty,y_1)$ and $J\cap(y_1,\infty)$, and using the right
part of (\ref{eq:trunc-hyp1}). Similarly, one has
$$
r_{\varphi,I_1}\{\varphi(y_1)\} = y_1 - \inf I_1 + \lambda\{t\in
I_1\cap(y_1,\infty) :\varphi(t)> \varphi(y_1)\},
$$
so that the following equality holds:
\begin{eqnarray}\label{eq:ry1}
r_{\varphi,J}\{\varphi(y_1)\} + \inf J & = &
r_{\varphi,I_1}\{\varphi(y_1)\}  + \inf I_1.
\end{eqnarray}

Now, define
$$
y_\star:= r_{\varphi,J}\{\varphi(y_1)\} + \inf J =
r_{\varphi,I_1}\{\varphi(y_1)\}  + \inf I_1.
$$
 It follows from this
definition that $T_J(\varphi)(y_\star) = \varphi(y_1) =
T_{I_1}(\varphi)(y_\star)$. Besides, $y_\star \leq \inf I_0$. To
prove this, note that  $T_J(\varphi)(\inf I_0) \leq M$  because of
the right parts of (\ref{eq:trunc-hyp1}) and (\ref{eq:trunc-hyp2});
so $y_\star \leq \inf I_0$ will follow as soon as $
T_J(\varphi)(y_\star) \geq M$, since $ T_J(\varphi)$ is a decreasing
function. This last inequality can be proved easily by
contradiction, using jointly that $T_J(\varphi)(y_\star) =
\varphi(y_1)$ and the left part of (\ref{eq:trunc-hyp2}).

 Finally, let us check that if both functions $T_J(\varphi)$ and
$T_{I_1}(\varphi)$ cross at one point (say, $y_\star$), then they
will coincide for each point $ \sup I_0 \geq x \geq y_\star $: Under
the hypothesis that they cross at  $y_\star$, it is equivalent to
show that for each $ -M \leq u \leq \varphi(y_1)$, one gets
\begin{eqnarray}\label{eq:ru}
 r_{\varphi,J}(u) + \inf J  =  r_{\varphi,I_1}(u) + \inf I_1.
 \end{eqnarray}
Let $u \in  [-M, \varphi(y_1)]$, and write on one hand
\begin{eqnarray*}
r_{\varphi,J}(u)& = & \lambda\{t\in J:\varphi(t)> u\}\\
& =& \lambda\{t\in J :\varphi(t)> \varphi(y_1)\} +
\lambda\{t\in J :\varphi(y_1) \geq \varphi(t)> u \}\\
&=& r_{\varphi,J}\{\varphi(y_1)\} + \lambda\{t\in J\cap(y_1,\infty)
:\varphi(y_1) \geq \varphi(t)> u \}\\
& = & r_{\varphi,I_1}\{\varphi(y_1)\} + \inf I_1 - \inf J\\
&& +
\lambda\{t\in I_1\cap(y_1,\infty) :\varphi(y_1) \geq \varphi(t)> u
\},
\end{eqnarray*}
where the last equality follows from (\ref{eq:ry1}) and the right
part of (\ref{eq:trunc-hyp1}). On the other hand,
\begin{eqnarray*}
r_{\varphi,I_1}(u) & = & r_{\varphi,I_1}\{\varphi(y_1)\} +
\lambda\{t\in
I_1 :\varphi(y_1) \geq \varphi(t)> u \}\\
& = & r_{\varphi,I_1}\{\varphi(y_1)\} + \lambda\{t\in
I_1\cap(y_1,\infty) :\varphi(y_1) \geq \varphi(t)> u \},
\end{eqnarray*}
so that equality (\ref{eq:ru}) holds, and this concludes the proof
of Theorem \ref{thm:trunc}.
 \eop Theorem \ref{thm:trunc} implies  that an extension of the definition of
$T_I$ to $I={\mathbb R}$ can be given for any continuous function
$\varphi\in {\cal B}({\mathbb R})$ such that (\ref{eq:trunc-hyp1})
and (\ref{eq:trunc-hyp2}) hold. Indeed, for any finite interval $J$ big
enough, $T_J(\varphi)(t)$ does not depend anymore on $J$, so that
one can define, for each $t \in I_0$:
$$
T(\varphi)(t):= T_{I_1}(\varphi)(t).
$$
A straightforward consequence of this definition is that both Lemma
\ref{lem:proprea}, Theorem \ref{thm:contraction} and Theorem
\ref{thm:trunc} hold for $T$:

\bcor Let $I_0\subset{\mathbb R}$ be a finite and fixed interval.
Assume $\varphi$ is continuous and satisfies (\ref{eq:trunc-hyp1})
and (\ref{eq:trunc-hyp2}). Then
\begin{eqnarray*}
(i)&&\mbox{$T$ satisfies Lemma \ref{lem:proprea} with the equalities and
inequalities assumed}\\
&&\mbox{to hold on $I_0$,}\\
(ii)&&\mbox{$T$ satisfies Theorem \ref{thm:contraction} with norm $||\cdot||$
defined on the set of}\\
&&\mbox{functions on $I_0$,}\\
(iii)&&\mbox{Theorem \ref{thm:trunc} holds with $T_{J}$ replaced by
$T$.}
\end{eqnarray*}
\ecor
\par
\section{The monotone estimation procedure}
\subsection{Definition and first properties}
Let $x$ be a function of interest (such as a density function, or a
regression function) and assume $x$ is non increasing. Consider an
estimator $x_n$ of $x$ constructed from $n$ observations, which is
not supposed to be monotone. Typically, $x_n$ can be an estimator
based on kernel, wavelets, splines, etc.
\bd We define as a new estimator of $x$ the monotone rearrangement
of $x_n$, namely $T(x_n)$. This is a non increasing estimator of
$x$. \ed

 \bth \label{thm:cv-estim} (i). Assume
that $\{x_n\}_{n\geq 1}$ is a uniformly consistent estimator of $x$
(in probability, uniformly on a compact set $B\subset {\mathbb R}$).
If $x$ is non increasing, then $\{T(x_n)\}_{n\geq 1}$ is a uniformly
consistent estimator of $x$ (in probability, uniformly on $B$).
\par
(ii). Assume that $\{x_n\}_{n\geq 1}$ is an estimator that converges
in probability in ${\mathbb L}^p$ norm to $x$. If $x$ is non
increasing, then $\{T(x_n)\}_{n\geq 1}$ converges in probability in
${\mathbb L}^p$ norm to $x$. \eth \prf  Both $(i)$ and $(ii)$ follow
from the fact that $||x||=\sup_{t\in K}|x(t)|$ is a norm, and $T$ a
contraction with respect to $||\cdot||$, by Theorem 1. Moreover
$T(x)=x$ if $x$ is non increasing. \eop

 \bnote  The
strong convergence in ${\mathbb L}^p$-norm of $T(f_n)$ to $f$, as a
consequence of the corresponding result for $f_n$, was first
established in Foug\`eres \cite[Theorem 5]{fougeres:1997} in the
case when $f_n$ is the kernel estimator of a density function $f$.
Chernozhukov et al. \cite{chernozhukov:fernandezval:galichon:2007a}
give a refinement of the non expansivity property, see their
Proposition 1, part 2,  providing a bound for the gain done by
rearranging $f_n$ and examining the multivariate framework as well.
\enote

\par
\subsection{Limit distribution results}
 Let $J\subset {\mathbb R}$ be a finite or infinite interval, and $C(J)$ the set of continuous functions on $J$.
 Let $x_n$ be a stochastic process in $C(J)$ and let $t_0$ be a fixed interior point in $J$. In this section
 limit distribution results for the random variable $T(x_n)(t_0)$ will be derived, where  $T$ is the
 monotone rearrangement map. The proof of these results are along the lines of Anevski and
  H\"ossjer \cite{anevski:hossjer:2006},
 and their notation will be used for clarity.
\par
Assume that $\{x_{n}\}_{n\geq 1}$ is a sequence of stochastic processes in $C(J)$ and write
\begin{eqnarray}
        x_{n}(t)&=&x_{b,n}(t)+v_{n}(t),\label{eq:PaAPROC}
\end{eqnarray}
for $t\in J$.
Given a sequence $d_{n}\downarrow 0$ and an interior point $t_{0}$ in $J$ define $ J_{n,t_{0}}=d_{n}^{-1}(J-t_{0})$.
Then, for $s\in J_{n,t_{0}}$, it is  possible to rescale the deterministic and stochastic parts of $x_n$ as
\begin{eqnarray*}
    \tilde{w}_{n}(s;t_{0})&=&d_{n}^{-1}\{v_{n}(t_{0}+sd_{n})-v_{n}(t_{0})\},\\
\tilde{g}_{n}(s)&=&d_{n}^{-1}\{x_{b,n}(t_{0}+sd_{n})-x_{b,n}(t_{0})\}.\label{eq:PaA4:1}
\end{eqnarray*}
which decomposes the rescaling of $x_n$ as
\begin{eqnarray*}
   d_n^{-1}\left\{x_n(t_0+sd_n)-x_n(t_0)\right\}&=&\tilde{g}_{n}(s)+\tilde{w}_{n}(s;t_{0}).
\end{eqnarray*}
However, due to the fact that the final estimator needs to be centered at the estimand $x(t_0)$ and not
at the preliminary estimator $x_n(t_0)$, it is more convenient to introduce the following rescaling
\begin{eqnarray}
    \tilde{v}_{n}(s;t_{0})&=&d_{n}^{-1}v_{n}(t_{0}+sd_{n}) \\
    &=&\tilde{w}_{n}(s;t_{0})+d_n^{-1}v_n(t_0),\nonumber\\
{g}_{n}(s)&=&d_{n}^{-1}\{x_{b,n}(t_{0}+sd_{n})-x(t_{0})\}\label{eq:PaA4:1}\\
&=&\tilde{g}_n(s)+d_n^{-1}\{x_{b,n}(t_0)-x(t_0)\},\nonumber
\end{eqnarray}
so that
\begin{eqnarray}
   y_n(s)&:=&g_n(s)+\tilde{v}_n(s; t_0)=d_n^{-1}\{x_{n}(t_0+sd_n)-x(t_0)\}.\label{eq:PaAyn}
\end{eqnarray}
This definition of the rescaled deterministic and stochastic parts
is slightly different from the one in Anevski and H\"ossjer
\cite{anevski:hossjer:2006}, and is due to the fact that we only
treat the case when the preliminary estimator and the final
estimator have the same rates of convergence, in which case our
definition is more convenient, whereas in Anevski and H\"ossjer
\cite{anevski:hossjer:2006} other possibilities occur.
\par
The limit distribution results will be derived using a classical
two-step procedure, cf. e.g. Prakasa Rao \cite{prakasarao:1969} :  A
local limit distribution is first obtained, under Assumption 1,
stating that the estimator $T(x_n)$ converges weakly in a local and
shrinking neighbourhood around a fixed point. Then it is shown,
under Assumption 2, that the limit distribution of $T(x_n)$ is
entirely determined by its behaviour in this shrinking
neighbourhood.
 \par
\bass There exists a stochastic process $\tilde{v}(\cdot;t_{0})\neq 0$ such that
\begin{eqnarray*}
    \tilde{v}_{n}(\cdot ;t_{0})&\stackrel{\cal L}{\rightarrow}& \tilde{v}(\cdot ;t_{0}),
\end{eqnarray*}
on $C(-\infty,\infty)$ as $n\rightarrow\infty$. The functions
$\{x_{b,n}\}_{n\geq 1}$ are monotone and there are constants $A<0$
and $\Delta \in{\mathbb R}$ such that for each $c>0$,
\begin{eqnarray}
\sup_{|s|\leq
c}|g_{n}(s)-(As+\Delta)|&\rightarrow&0,\label{eq:PaAbiasconv}
\end{eqnarray}
as $n\rightarrow\infty$.
\eass
In the applications typically
\begin{eqnarray*}
   A&=&\lim_{n\rightarrow \infty} \frac{\tilde{g}_n(s)}{s}=x'(t_0),\\
   \Delta &=&\lim_{n\rightarrow \infty} d_n^{-1}\{x_{b,n}(t_0)-x(t_0)\},
 \end{eqnarray*}
 so that $A$ is the local asymptotic  linear term and $\Delta$ is the local asymptotic bias, both
 properly normalized, of the preliminary estimator $x_{n}$.
 \par
 Define the (limit) function
\begin{eqnarray}
 y(s)&=&As+\Delta +\tilde{v}(s;t_{0}).\label{eq:PaAy}
\end{eqnarray}
Let $\{z_n\}$ be an arbitrary sequence of stochastic processes.
%
% ---- Assumption 2  --------
%
\bass  Let $I_0$ be a given compact interval and $\delta>0$.  There
exists a positive constant $c$ such that $[-c,c]\supset I_0$ and a
finite positive $M$ such that
\begin{eqnarray}
\liminf_{n\rightarrow\infty}P\left\{\inf_{s \in (-c,\sup
I_0)}z_{n}(s) \geq -M,\sup_{s\in (c,\infty)}z_n(s) \leq -M\right\}&>
&1-\delta,  \label{eq:downdip1}
\end{eqnarray}
and
\begin{eqnarray}
   \liminf_{n\rightarrow \infty} P\left\{\inf_{s\in (-\infty,-c)} z_{n}(s) \geq +M,
   \sup_{s\in (\inf I_0,c)} z_n(s)\leq +M \right\}&>&1-\delta. \label{eq:downdip2}
\end{eqnarray}

\eass
Denote $T_{c}=T_{[-c,c]}$ and  $T_{c,n}=T_{[t_{0}-cd_{n},t_{0}+cd_{n}]}$. The truncation
result Theorem 2 has a probabilistic counterpart in
the following.
 \blem\label{lem:A2} Let
$\{z_{n}\}$ satisfy Assumption 2. Let $I$ be a finite interval in
${\mathbb R}$. Then for each compact interval $J \supset [-c,c]$
\begin{eqnarray}
\lim_{c\rightarrow\infty}\limsup_{n\rightarrow\infty} P(\sup
_{I}|T_{c}(z_{n})(\cdot)-T_J(z_{n})(\cdot)|=0)&=&1.\nonumber
\end{eqnarray}
\elem
\prf Let $A_n$ and $B_n$ be the sets for which the probabilities are
bounded in $(\ref{eq:downdip1})$ and $(\ref{eq:downdip2})$,
respectively. Then, using Theorem $\ref{thm:trunc}$ with
$I_1=[-c,c]$ and $I_0=I$, it follows that $ A_n\cap B_n\subset \{
\sup_{I}| T_c(z_n)-T_J(z_n)|=0\}$ for each compact interval $J
\supset [-c,c]$. Since $(\ref{eq:downdip1})$ and
$(\ref{eq:downdip2})$ imply $P(A_n\cap B_n)\geq 1-2\delta$ it
follows that $\limsup_{n\rightarrow\infty} P(\sup
_{I}|T_{c}(z_{n})(\cdot)-T_J(z_{n})(\cdot)|=0)\geq 1-2\delta$. Since
$\delta>0$ is arbitrary, taking the limit as $c\to \infty$ of the
left hand side of this expression implies the statement of the
lemma.
 \eop
  Note that
the previous lemma holds with $T_J$ replaced by $T_{I_n}$ for an
arbitrary sequence of intervals $I_n$ growing to ${\mathbb R}$.
\bth\label{thm:asympt-distr} Let $J \subset {\mathbb R}$ be an
interval, and $t_{0}$ be a fixed point belonging to the interior of
$J$. Suppose Assumption 1 holds. Assume moreover that Assumption 2
holds for both $\{y_n\}$ and $y$. Then
\begin{eqnarray} {\label{eq:PaAj}}
d_{n}^{-1}\Large{[}T_{J}(x_{n})(t_{0})-x(t_{0})\Large {]}
&\stackrel{\cal L}{\rightarrow}& T  \Large{[} A\cdot+ \tilde{v}
    (\, \cdot \, ;t_{0})\Large{]}(0)+\Delta ,
\end{eqnarray}
as $n\rightarrow\infty$. \eth
\prf   Let $c>0$ be fixed. We have
\begin{eqnarray}
d_{n}^{-1}\{T_{J}(x_{n})(t_{0})-x(t_{0})\}&=&
d_{n}^{-1}\{T_{J}(x_{n})(t_{0})-T_{c,n}(x_{n})(t_{0})\}\nonumber
\\&&+d_{n}^{-1}\{T_{c,n}(x_{n})(t_{0})-x(t_{0})\}.
\label{eq:display01}
\end{eqnarray}
\par
Let us first consider the second term of the right hand side of $(\ref{eq:display01})$ and introduce
\begin{eqnarray}
  \chi_n(s)&:=&x_{n}(t_{0}+sd_{n})=x(t_0)+d_{n}y_n(s). \label{eq:PaAi}
\end{eqnarray}
Applying Lemma \ref{lem:proprea} $(i)$ and $(ii)$ leads to
\begin{eqnarray*}
     T_{c,n}(x_n)(t_0+sd_n)&=&T_c(\chi_n)(s)=d_nT_{c}(y_n)(s)+x(t_0),
\end{eqnarray*}
which gives
\begin{eqnarray}
d_{n}^{-1}\{T_{c,n}(x_{n})(t_{0})-x(t_{0})\}&= & T_{c}(y_n)(0).
\nonumber
\end{eqnarray}
Assumption 1 implies that $y_n\stackrel{\cal L}{\rightarrow} y$ on
$C[-c,c]$, with $y$ defined in $(\ref{eq:PaAy})$.  Applying the
continuous mapping theorem on $T_c$, cf. Theorem \ref{thm:contraction}, proves
\begin{eqnarray*}
d_Ä{n}^{-1}\{T_{c,n}(x_{n})(t_{0})-x(t_{0})\}&\stackrel{\cal
L}{\rightarrow}&T_{c} (y)(0)
\end{eqnarray*}
as $n\rightarrow\infty$. Lemma  \ref{lem:A2} via Assumption 2 with
$z_{n}=y$ implies
\begin{eqnarray*}
 T_{c}(y)(0)-T(y)(0)&\stackrel{P}{\rightarrow}&0
\end{eqnarray*}
as $c\rightarrow\infty$.
\par
Next we consider the first  term of the right hand side of $(\ref{eq:display01})$. Let $\nabla$ be a positive and
finite constant and denote $A_{n,\nabla}=[t_0-\nabla d_n,t_0+\nabla
d_n]$. From $(\ref{eq:PaAi})$ and Lemma \ref{lem:proprea} $(i,ii)$
it follows that
\begin{eqnarray*}
     \sup_{A_{n,\nabla}}d_{n}^{-1}|T_{c,n}(x_{n})(\cdot)-T_{J}(x_{n})(\cdot)|
     = \sup_{[-\nabla,\nabla]}|T_{c}(y_{n})(\cdot)-T_{J_{n,t_{0}}}(y_{n})(\cdot)|,
\end{eqnarray*}
with $y_{n}$ as defined in $(\ref{eq:PaAyn})$. If $J= {\mathbb R}$
(resp. $J\ne {\mathbb R}$),  Lemma \ref{lem:A2} (resp. note
following Lemma \ref{lem:A2}) can be used with $I=[-\nabla,\nabla]$
to obtain
\begin{eqnarray*}
d_{n}^{-1}\{T_{c,n}(x_{n})(t_{0})-T_{J}(x_{n})(t_{0})\}&\stackrel{P}{\rightarrow}&0
\end{eqnarray*}
 if we first let $n\rightarrow\infty$ and then let $c\rightarrow\infty$.
\par
Letting first $n$ and then $c$ tend to infinity in $(\ref{eq:display01})$,
 applying Slutsky's theorem and Lemma 3 $(i)$ finishes the proof.
\eop
 \bnote
The approach for deriving the limit distributions is similar to the
general approach in Anevski and H\"ossjer
\cite{anevski:hossjer:2006} with a preliminary estimator that is
made monotone via the ${\mathbb L}^2$-projection on the space of
monotone functions. There are however a few differences: \\
$-$
Anevski and H\"ossjer look at rescaling of an integrated preliminary
estimator of the monotone functions, whereas we rescale the
estimator directly. Our approach puts a stronger assumption on the
asymptotic properties of the preliminary estimator, which is however
traded off against weaker conditions on the map $T$, since we
only have to assume that the map $T$ is continuous; had we dealt
with rescaling as in Anevski and H\"ossjer we would have had to
prove that the composition $ \frac{d}{dt}(\tilde{T})$ (with
$\tilde{T}$ defined by $\tilde{T}(F)(t)=\int^t_0 T(F')(u)\,du$) is a
continuous map, which is generally not true for $T$ equal to the
monotone rearrangement map; it is however true, under certain
conditions, for $\tilde{T}$ equal to the least concave minorant map
(when $T$ becomes the ${\mathbb L}^2$-projection on the space of
monotone functions), cf. Proposition 2 in Anevski and H\"ossjer
\cite{anevski:hossjer:2006}.\\
$-$ Furthermore, we are able to do rescaling for the preliminary
estimator directly since it is a smooth function. On the contrary,
for some of the cases treated in Anevski and H\"ossjer this is not
possible, e.g. for the isotonic regression and the NPMLE of a
monotone density the rescaled stochastic part is asymptotically
white noise. As a consequence our rescaled deterministic function is
assumed to be approximated by a linear function, whereas the
rescaled deterministic function in Anevski and H\"ossjer
\cite{anevski:hossjer:2006} is assumed to be
approximated by a convex or concave function. \\
$-$ Finally, the
rescaling is here centered at $x(t_0)$, and not at $x_n(t_0)$, which
makes it more convenient to apply the limit distribution result we
get.
 \eop \enote

\section{Applications to nonparametric inference problems}
In this section we present estimators of a monotone density function
 and a monotone regression function. Limit distributions for
estimators of a marginal decreasing density $f$ for stationary
weakly dependent data with marginal  density $f$ as well as of a
monotone regression function $m$ with stationary errors, that are
weakly or strongly dependent, will be derived.
\par
For the density estimation problem let  $\{t_i\}_{i=1}^{\infty}$
denote a stationary process with marginal density function $f$.
Define the empirical distribution function
$F_{n}(t)=\frac{1}{n}\sum_{i=1}^{n} 1_{\{t_{i}\leq t\}}$ and the
centered empirical process
$F_{n}^{0}(t)=\frac{1}{n}\sum_{i=1}^{n}(1_{\{t_{i}\leq t\}}-F(t))$.
Consider a sequence $\delta_{n}$ such that $\delta_{n}\downarrow
0,n\delta_{n}\uparrow \infty$ as $n\rightarrow\infty$, and define
the centered empirical process locally around $t_{0}$ on scale
$\delta_{n}$ as
\begin{eqnarray}
w_{n,\delta_{n}}(s;t_{0})&=&\sigma_{n,\delta_{n}}^{-1}n\{F_{n}^{0}(t_{0}+s\delta_{n})-F_{n}^{0}(t_{0})\}\nonumber\\
&=&\sigma_{n,\delta_{n}}^{-1}\sum_{i=1}^{n}(1_{\{t_{i} \leq
t_{0}+s\delta_{n}\}}-1_{\{t_{i}\leq t_{0}\}}
\nonumber\\
&&-F(t_{0}+s\delta_{n})+F(t_{0})),\label{eq:cent_emp_proc}
\end{eqnarray}
where
\begin{eqnarray*}
  \sigma_{n,\delta_{n}}^{2}&=&\mbox{Var}\left[n\left\{F_{n}^{0}(t_{0}+\delta_{n})-F_{n}^{0}(t_{0})\right\}\right]\\
   &=&\mbox{Var}\left[\sum_{i=1}^{n}\left\{1_{\{t_{0} < t_{i}\leq t_{0}+\delta_{n}\}}
   -F(t_{0}+\delta_{n})+F(t_{0})\right\}\right].
\end{eqnarray*}
\par
For the regression function estimation problem let
$\{\epsilon_{i}\}_{i=-\infty}^{\infty}$ be a stationary sequence of
random variables with $E(\epsilon_{i})=0$ and
Var$(\epsilon_{i})=\sigma^{2}<\infty$. Let
$\sigma_{n}^{2}=$Var$(\sum_{i=1}^{n}\epsilon_{i})$. The two sided
partial sum process $w_{n}$ is defined by
\begin{eqnarray}
           w_{n}(t_{i}+\frac{1}{2n})&=&\left\{ \begin{array}{ll}
                                  \frac{1}{\sigma_{n}}(\frac{\epsilon_{0}}{2}+\sum_{j=1}^{i}\epsilon_{j}),\;\;\;\;i=0,1,2,\ldots,\nonumber\\
                                  \frac{1}{\sigma_{n}}(-\frac{\epsilon_{0}}{2}-\sum_{j=i+1}^{-1}\epsilon_{j}),\;\;\;\;i=-1,-2,\ldots,
                                         \end{array} \right .
\end{eqnarray}
and linearly interpolated between these points. Note that $w_n\in C({\mathbb R})$.
\par
Let $\Cov(k)=E(\xi_{1}\xi_{1+k})$  denote the covariance function of
a generic stationary sequence $\{\xi_i\}$, and distinguish between
three cases (of which {\bf [a]} is a special case of {\bf [b]}.)
\begin{enumerate}
 \item[{\bf [a]}] Independence: the $\epsilon_{i}$ are independent.
 \item[{\bf [b]}] Weak dependence: $\sum_{k}|\Cov(k)|<\infty$.
 \item[{\bf [c]}] Strong (long range) dependence: $\sum_{k}|\Cov(k)|=\infty$.
\end{enumerate}
\par
Weak dependence can be further formalized using mixing conditions as
follows: Define two $\sigma$-algebras of a sequence $\{\xi_i\}$ as
${\cal
  F}_{k}=\sigma\{\xi_{i}: i\leq k\}$ and
$ \bar{{\cal  F}}_{k}=\sigma\{\xi_{i}:i\geq k\}$, where
$\sigma\{\xi_{i}:i\in I\}$ denotes the $\sigma-$algebra generated by
$\{\xi_{i}:i\in I\}$. The stationary sequence $\{\xi_{i}\}$ is said
to be "$\phi$-mixing" or "$\alpha$-mixing" respectively if there is
a function $\phi(n)$ or $\alpha(n)\rightarrow 0$ as $n\rightarrow
\infty$, such that
\begin{eqnarray}\label{eq:mixing}
   \sup_{A\in \bar{{\cal F}}_{n}}|P(A|{\cal
  F}_{0})-P(A)|&\leq&\phi(n), \nonumber \\
   \sup_{A\in {\cal F}_{0},B\in \bar{{\cal F}}_{n}}|P(AB)-P(A)P(B)|&\leq
   &\alpha(n),
\end{eqnarray}
respectively.
\par
Long range dependence is usually formalized using subordination or assuming the processes are linear;
we will treat only (Gaussian) subordination.
\par
All limit distribution results stated will be for processes in $C(-\infty,\infty)$ with the uniform
metric on compact intervals and the Borel $\sigma$-algebra.

 \subsection{Monotone regression function estimation}
 In this section we introduce an estimator of  a monotone regression function. We derive
  consistency and limit distributions, under general dependence assumptions.
 \par
Assume $m$ is a $C^1$-function on a compact interval $J\subset
{\mathbb R}$, say $J=[0,1]$ for simplicity; let
$(y_{i},t_{i}),i=1,\cdots,n$ be pairs of data satisfying
\begin{eqnarray}\label{eq:PaAregr}
y_{i}&=&m(t_{i})+\epsilon_{i},
\end{eqnarray}
where $t_{i}=i/n$.
\par
Define $\bar{y}_{n}:[1/n,1]\mapsto {\mathbb R}$ by linear
interpolation of the points $\{(t_{i},y_{i})\}_{i=1}^{n}$, and let
\begin{eqnarray}\label{eq:PaAkernel}
x_{n}(t)&=&h^{-1} \int k((t-u)/h)\bar{y}_{n}(u)\,du,
\end{eqnarray}
be the Gasser-M\"uller kernel estimate of $m(t)$, cf. Gasser and
M\"uller \cite{gasser:muller:1984}, where $k$ is a density in ${\bf
L}^{2}({\mathbb R})$ with compact support, for simplicity take
supp$(k)=[-1,1]$. Let $h$  be the bandwidth, for which we assume
that $h\rightarrow 0,nh\rightarrow\infty$.
\par
 To define a monotone estimator of $m$, we put
\begin{eqnarray}
        \tilde{m}(t)&=&T_{[0,1]}(x_{n})(t),\;t\in J,\label{eq:tildem}
\end{eqnarray}
where $T$ is the monotone rearrangement map. A straightforward
application of Theorem \ref{thm:cv-estim} and standard consistency results for regression function estimators imply the
following consistency result:
\bpr The random function $\tilde{m}$ defined by (\ref{eq:tildem}) is
a uniformly consistent estimator of $m$ in probability uniformly on
compact sets, and in probability in ${\mathbb L}^p$ norm. \epr
 Clearly $x_{n}(t)=x_{b,n}(t)+v_{n}(t)$, with
\begin{eqnarray}
x_{b,n}(t)&=&h^{-1}\int k(\frac{t-u}{h})\bar{m}_{n}(u)\,du, \label{def:xbn-reg} \\
v_{n}(t)&=&h^{-1}\int k(\frac{t-u}{h})\bar{\epsilon}_{n}(u)\,du,
\nonumber
\end{eqnarray}
where the functions $\bar{m}_{n}$ and $\bar{\epsilon}_{n}$ are
obtained by linear interpolation of \linebreak $\{(t_{i},m(t_{i}))\}_{i=1}^{n}$
and $\{(t_{i},\epsilon_{i})\}_{i=1}^{n}$ respectively. For the
deterministic term \linebreak $x_{b,n}(t)\rightarrow x_{b}(t)=m(t)$, as
$n\rightarrow\infty$. Note that $\bar{m}_{n}$, and thus also
$x_{b,n}$, is monotone.
\par
Put
\begin{eqnarray}
      \bar{w}_{n}(t)&=&\frac{n}{\sigma_{n}}\int_{0}^{t}\bar{\epsilon}_{n}(u)\,du.\label{eq:PaAwntilde}
\end{eqnarray}
Since supp$(k)=[-1,1]$ and if $t\in(1/n+h,1-h)$, from a partial integration and change of variable we obtain
\begin{eqnarray*}
      v_{n}(t)&=&\frac{\sigma_{n}}{nh}\int k'(u)\bar{w}_{n}(t-uh)\,du.
\end{eqnarray*}
It can be shown that $\bar{w}_{n}$ and $w_{n}$ are asymptotically
equivalent for all dependence structures treated in this paper. Let
us now recall how the two sided partial sum process behaves in the
different cases of dependence we consider:
%\begin{itemize}
%\item[

{\bf [a]} When the $\epsilon_{i}$ are independent, we have the
classical Donsker theorem, cf. $\!$Billingsley
\cite{billingsley:1968}, implying that
\begin{eqnarray}
        w_{n}\stackrel{\cal L}{\rightarrow} B, \label{eq:PaAweak1}
\end{eqnarray}
as $n\rightarrow\infty$, with $B$ a two sided standard Brownian
motion on $C({\mathbb R})$.
%\item[

{\bf [b]} Define
\begin{eqnarray}
  \kappa^{2}&=&\Cov(0)+2\sum_{k=1}^{\infty}\Cov(k). \label{eq:PaAkappa2}
\end{eqnarray}
\bass $[{\bf \phi-mixing}]$ Assume $\{\epsilon_{i}\}_{i\in {\bf Z}}$
is a stationary $\phi$-mixing sequence with $E\epsilon_{i}=0$ and
$E\epsilon_{i}^{2}<\infty$. Assume further $\sum_{k=1}^{\infty}
\phi(k)^{1/2}<\infty$ and $\kappa^2>0$ in $(\ref{eq:PaAkappa2})$.
\eass
\bass $[{\bf \alpha-mixing}]$ Assume $\{\epsilon_{i}\}_{i\in {\bf
Z}}$ is a stationary $\alpha$-mixing sequence with $E\epsilon_{i}=0$
and  $E\epsilon_{i}^{4}<\infty$, $\kappa^2>0$ in
$(\ref{eq:PaAkappa2})$ and  $\sum_{k=1}^{\infty}
\alpha(k)^{1/2-\epsilon}<\infty$, for some $\epsilon>0$. \eass
Assumption 3 or 4 imply that $\sigma_n^2\rightarrow \kappa^2$ and
that Donsker's result $(\ref{eq:PaAweak1})$ is valid, cf. Anevski
and H\"ossjer \cite{anevski:hossjer:2006} and references therein.
%\item[

{\bf [c]} We model long range dependent data $\{\epsilon_i\}_{i\geq
1}$ using Gaussian subordination: More precisely, we write
$\epsilon_i=g(\xi_i)$ with $\{\xi_{i}\}_{i\in {\bf Z}}$ a stationary
Gaussian process with mean zero and covariance function
$\Cov(k)=E(\xi_{i}\xi_{i+k})$ such that $\Cov(0)=1$ and
$\Cov(k)=k^{-d}l_{0}(k)$, with $l_{0}$ a slowly varying function  at
infinity\footnote{i.e. $l_0(tk)/l_0(t) \to 1$ as $t\to \infty$ for
each positive $k$.} and $0<d<1$ fixed. Furthermore $g:{\mathbb
R}\mapsto {\mathbb R}$ is a measurable function with
$E\{g(\xi_1)^2\}<\infty$. An expansion $g(\xi_{i})$ in Hermite
polynomials is available
\begin{eqnarray*}
     g(\xi_{i})&=&\sum_{k=r}^{\infty}\frac{1}{k!}\eta_{k}h_{k}(\xi_{i}),
\end{eqnarray*}
where equality holds as a limit in  $L^{2}(\varphi)$, with $\varphi$
the standard Gaussian density function. The functions
$h_{k}(t)=t^{-k}(d/dt)^k(t^ke^{-t^2})$ are the Hermite polynomials
of order $k$, the functions
\begin{eqnarray*}
\eta_{k}=E\left\{g(\xi_{1})h_{k}(\xi_{1})\right\}=\int
g(u)h_{k}(u)\phi(u)\,du, \nonumber
\end{eqnarray*}
 are the $L^{2}(\varphi)$-projections on $h_{k}$, and $r$ is the index of the first non-zero
 coefficient in the expansion. Assuming that $0<dr<1$, the subordinated sequence
 $\{\epsilon_{i}\}_{i\geq 1}$ exhibits long range dependence (see e.g. Taqqu
 \cite{taqqu:1975,taqqu:1979}), and Taqqu \cite{taqqu:1975} also
 shows that
\begin{eqnarray}
  \sigma_{n}^{-1}\sum_{i\leq nt}g(\xi_{i})&\stackrel{\cal L}{\rightarrow}& z_{r
,\beta}(t), \nonumber
\end{eqnarray}
in $D[0,1]$ equipped with the Skorokhod topology, with variance
$\sigma_{n}^{2}=$ Var \linebreak
$\left\{\sum_{i=1}^{n}g(\xi_{i})\right\} =
\eta_{r}^{2}n^{2-rd}l_{1}(n)(1+o(1))$, where
\begin{eqnarray}\label{eq:PaAl1}
  l_{1}(k)&=&\frac{2}{r!(1-rd)(2-rd)}l_{0}(k)^{r}.
\end{eqnarray}
 The limit process $z_{r,\beta}$ is in $C[0,1]$ a.s., and is self similar with parameter
\begin{eqnarray}
\beta&=&1-rd/2. \label{eq:PaAbetadef}
\end{eqnarray}
The process $z_{1,\beta}(t)$ is fractional Brownian motion,
$z_{2,\beta}(t)$ is the Rosenblatt process, and the processes
$z_{r,\beta}(t)$ are all non-Gaussian for $r\geq 2$, cf. $\!$Taqqu
\cite{taqqu:1975}.  From these results follows a two sided version
of Taqqu's result stating the behavior of the two sided partial sum
process:
\begin{eqnarray}
   w_{n}&\stackrel{\cal L}{\rightarrow}& B_{r,\beta},\label{eq:PaAweak3}
\end{eqnarray}
in $D(-\infty,\infty)$, as $n\rightarrow\infty$, where $B_{r,\beta}$
are the two sided versions of the processes $z_{r,\beta}$.
%\end{itemize}
\par
In the sequel,  rescaling is done at the bandwidth rate, so that
$d_{n}=h$. For $s>0$, let consider the following  rescaled process:
\begin{eqnarray}
  \tilde{v}_{n}(s;t)&=&d_{n}^{-1}(nh)^{-1}\sigma_{\hat{n}}\int\bar{w}_{\hat{n}}(h^{-1}t+s-u) k'(u)\,du \nonumber \\
&\stackrel{\cal
L}{=}&d_{n}^{-1}(nh)^{-1}\sigma_{\hat{n}}\int\bar{w}_{\hat{n}}(s-u)
k'(u)\,du, \label{eq:vnrescaled}
\end{eqnarray}
with $\hat{n}=[nh]$ the integer part of $nh$, where the last
equality holds due to the stationarity (exactly only for $t=t_{i}$
and asymptotically otherwise). Note that the right hand side holds
also for $s<0$.

With the bandwidth choice $d_n=h$ we obtain a non-trivial limit
process $\tilde{v}$; choosing $d_n$ such that $d_n/h\to 0$ leads to
a limit ``process'' equal to a random variable and $d_n/h\to \infty$
to white noise. In the first case the limit distribution of $T(x_n)$
on the scale $d_n$ will be the constant $0$, while in the second
case it will (formally) be $T(m'(t_0)\cdot+\tilde{v}(\cdot))(0)$
which is not defined ($T$ can not be defined for generalized
functions, in the sense of L. Schwartz \cite{schwartz:1966}).

\bth\label{thm:reglim} Assume $m$ is monotone on $[0,1]$ and for
some open interval $I_{t_0}\ni t_0$, $m\in C^{1}(I_{t_0})$ and
$\sup_{t\in I_{t_0}} m'(t)<0$ with $t_{0}\in (0,1)$.  Let $x_{n}$ be
the kernel estimate of $m$ defined in $(\ref{eq:PaAkernel})$,
    with a non-negative and compactly supported kernel $k$ such that
$k'$ is bounded, and with bandwidth $h$ specified below. Suppose that one of the following conditions holds.
\begin{eqnarray*}
 {\bf [a]} &&\mbox{$\{\epsilon_{i}\}$ are independent and identically
 distributed, $E\epsilon_{i}=0$;} \\
 %{\mathbb E}\epsilon_i^3 <\infty$
 %{\tt [Remove this hypothesis ?]}}\\
&& \mbox{$\sigma^{2}=$Var$(\epsilon_{i})<\infty$,} \mbox{ and $h=an^{-1/3}$, for an arbitrary $a>0$},\\
        {\bf [b]}&&\mbox{Assumption 3 or 4 holds, $\sigma_{n}^{2}=\Var(\sum_{i=1}^{n}\epsilon_{i})$,}
% && \mbox{ $\;{\mathbb E}\epsilon_i^3 <\infty$ ;
%{\tt [Remove this hypothesis ?]} }\\
 \mbox{ $\kappa^{2}$~is~defined~} \mbox{in
 $(\ref{eq:PaAkappa2})$},\\
  && \mbox{and $h=an^{-1/3}$, with $a>0$ an arbitrary constant},\\
        {\bf [c]}&&\mbox{$\epsilon_{i}=g(\xi_{i})$ is a long range dependent subordinated Gaussian sequence
}\\
        &&\mbox{with parameters $d$ and $r$, $h=l_{2}(n;a)n^{-rd/(2+rd)}$} \mbox{ with $a>0$ and}\\
        &&\mbox{$n\mapsto l_{2}(n;a)$ is a slowly varying function  defined in the proof below.}
\end{eqnarray*}
Then, correspondingly, we obtain
\begin{eqnarray*}
   h^{-1}\{\tilde{m}(t_0)-m(t_{0})\}&\stackrel{\cal L}{\rightarrow}&
    T[m'(t_{0})\cdot +\tilde{v}(\cdot ;t_{0})](0)+m'(t_0)\int
   uk(u)\,du,
\end{eqnarray*}
as $n\rightarrow\infty$, where $\tilde{m}$ is defined in
(\ref{eq:tildem}),
\begin{eqnarray}{\label{eq:PaAproc}}
\tilde{v}(s;t)&=&c\int w(s-u) k'(u)\,du,
\end{eqnarray}
and respectively
\begin{eqnarray*}
{\bf [a]} &&w=B\; ; \; c=\sigma a^{-3/2},\\
{\bf [b]}&& w=B\; ; \;c=\kappa a^{-3/2},\\
{\bf [c]}&& w=B_{r,\beta}\; ; \;c=|\eta_{r}|a ~~({\mbox{where~}}
\beta {\mbox{~defined~in~}}
 (\ref{eq:PaAbetadef})).
\end{eqnarray*}
\eth
\prf Theorem \ref{thm:reglim} is an application of Theorem
\ref{thm:asympt-distr} in the context of monotone regression
function. Assume first that $d_{n}=h$ is such that
\begin{eqnarray}{\label{eq:PaAconst}}
d_{n}^{-1}(nh)^{-1}\sigma_{\hat{n}}&=&d_{n}^{-2}n^{-1}\sigma_{\hat{n}}\rightarrow c>0.
\end{eqnarray}
Then  $w_{n}\stackrel{\cal L}{\rightarrow} w$ in
$D(-\infty,\infty)$, using the supnorm over compact intervals
metric, under the respective assumptions in ${\bf [a]},{\bf [b]}$
and ${\bf [c]}$. Besides, note that if $k'$ is bounded and $k$ has
compact support, the map
\begin{eqnarray*}
 C(-\infty,\infty)  \ni z(s)&\mapsto&\int z(s-u)k'(u)\,du \in C(-\infty,\infty)
\end{eqnarray*}
is continuous, in the supnorm over compact intervals metric. Thus,
under the assumptions that $k'$ is bounded and $k$ has compact
support,  the continuous mapping theorem implies that
\begin{eqnarray}{\label{eq:PaAproc-conv}}
\tilde{v}_{n}(s;t)&\stackrel{\cal L}{\rightarrow}&\tilde{v}(s;t),
\end{eqnarray}
where $\tilde{v}(s;t)$ is defined in $(\ref{eq:PaAproc})$. This
yields the first part of Assumption 1. Furthermore

\begin{eqnarray*}
        \tilde{g}_{n}(s)&=&h^{-1}\int \ell(u)\bar{m}_{n}(t_{0}-hu)\,du\\
        &=&h^{-1}\int \ell(u){m}(t_{0}-hu)\,du+r_{n}(s),
\end{eqnarray*}
with $\ell(v)=k(v+s)-k(v)$ and $r_n$ a remainder term. Since
\begin{eqnarray*}
        \int v^{\lambda}\ell(v)\,dv&=&\left\{ \begin{array}{ll}
                               0 ,\;  {\rm ~if~} \lambda=0, \\
                               -s,\; {\rm ~if~} \lambda=1,
                                       \end{array} \right.\\
\end{eqnarray*}
it follows by a Taylor expansion of $m$ around $t_{0}$ that the
first term converges towards $As$, with $A=m'(t_{0})$.  The
remainder term is bounded for any $c>0$ as
\begin{eqnarray*}
         \sup_{|s|\leq c}|r_{n}(s)|&\leq &h^{-1} \sup_{|s|\leq c} \int |\ell(u)|\,du\,
         \sup_{|u-t_{0}|\leq (c+1)h}|\bar{m}_{n}(u)-m(u)|\\
         &=&O(n^{-1}h^{-1})=o(1).
\end{eqnarray*}

 Furthermore
\begin{eqnarray}
       d_n^{-1}\{x_{b,n}(t_0)-m(t_0)\}&\rightarrow& m'(t_0)\int
   uk(u)\,du \label{eq:rev3_003}=:\Delta,
\end{eqnarray}
as $n\rightarrow\infty $, which proves Assumption 1.
\par
Proof that  Assumption 2 holds is relegated to the appendix, see
Corollary \ref{cor:reg} in Appendix \ref{append-reg}. An application
of Theorem~\ref{thm:asympt-distr} then finishes the proof of
Theorem~\ref{thm:reglim}. It only remains to check whether
$d_{n}^{-1}(nh)^{-1}\sigma_{\hat{n}}\to c>0$ for the three types of
dependence.
\par
\begin{itemize}
\item[-] Independent case ${\bf [a]}$: We have $\sigma^{2}_{\hat{n}}=\sigma^{2}n
d_{n}$. Thus $d_{n}^{-1}(nh)^{-1}\sigma_{\hat{n}}=\sigma
n^{-1/2}h^{-3/2}$, and $(\ref{eq:PaAconst})$ is satisfied with
$c=\sigma a^{-3/2}$ if $d_{n}=h=an^{-1/3}$.
\item[-] Mixing case ${\bf [b]}$: The proof is similar to the proof of ${\bf [a]}$,
replacing $\sigma$ by $\kappa$.
\item[-] Long range data case ${\bf [c]}$:  Since
$\sigma^{2}_{\hat{n}}=\eta_{r}^{2}(nd_{n})^{2-rd}l_{1}(nd_{n})$, if we choose  $d_{n}=h$  we will have
\begin{eqnarray}
  d_n^{-2}n^{-1}\sigma_{\hat{n}}=d_{n}^{-2}n^{-1}|\eta_{r}|(nd_{n})^{1-rd/2}l_{1}(nd_{n})^{1/2}
  &\to& |\eta_{r}|a \label{eq:PaA000}
\end{eqnarray}
if and only if
\begin{eqnarray}
     d_{n}&=&n^{-rd/(2+rd)}l_{2}(n;a)\label{eq:PaA0},
\end{eqnarray}
where $l_{2}$ is another function slowly varying at infinity,
implicitly defined in $(\ref{eq:PaA000})$. Thus
($\ref{eq:PaAconst}$) follows with $c=|\eta_{r}|a$ and $h=d_n$ given
in $(\ref{eq:PaA0})$.\hskip 2cm $\Box$
\end{itemize}
%\vskip -14mm \eop
%
%
\bnote The present estimator is similar to the estimator first
presented by Mammen \cite{mammen:1991}: Mammen proposed to do
isotonic regression  of a kernel estimator of a regression function (using bandwidth $h=n^{-1/5}$),
whereas we do monotone rearrangement of a kernel estimator. Mammen's estimator was extended to dependent data and other bandwidth choices by Anevski
and H\"ossjer \cite{anevski:hossjer:2006} who derived limit
distributions for weak dependent and long range dependent data that
are analogous to our results; for the independent data case and bandwidth choice $h=n^{-1/3}$ the
limit distributions are similar with rate of convergence $n^{1/3}$
and nonlinear maps of Gaussian processes.  \enote
\subsection{Monotone density estimation}
In this subsection we introduce a (monotone) estimator of a monotone
density function for stationary data, for which we derive
consistency and limit distributions.
\par
Let $t_1,t_2,\ldots$ denote a stationary process with marginal
density function $f$ lying in the class  of decreasing density
functions on ${\mathbb R^+}$, and define the following estimator
 of the marginal decreasing density for the
sequence $\{t_i\}_{i\geq 1}$: Consider
$x_{n}(t)=(nh)^{-1}\sum_{i=1}^{n}k\{({t-t_i})/{h}\}$  the kernel
estimator of the density $f$, with $k$ a bounded density function
supported on $[-1,1]$ such that $\int k'(u) du =0$, and $h>0$ the
bandwidth (cf. e.g. Wand and Jones \cite{wand:jones:1995}), and
define the (monotone) density estimate
\begin{eqnarray}
\hat{f}_n(t)&=&T(x_n)(t), \label{def:fnrea}
\end{eqnarray}
where $T$ is the monotone rearrangement map. Note that
$\hat{f}_n$ is monotone and positive, and integrates to one, cf.
equation (4) of Section 3.3. in Lieb and Loss \cite{lieb:loss:1996}.
\par
A straightforward consequence of Theorem \ref{thm:cv-estim} and
standard convergence results for the kernel density estimate is the following consistency result:
\bpr The random function $\hat{f}_n$ defined by (\ref{def:fnrea}) is
a uniformly consistent estimator of $f$ in probability uniformly on
compact sets, and in probability in ${\mathbb L}^p$ norm. \epr
In the following, the limit distributions for $\hat{f}_n$ in the
independent and weakly dependent cases are derived. We will in
particular  make use of recent results on the weak convergence
$w_{n,\delta_{n}}\stackrel{\cal L}{\rightarrow}w$, on
$D(-\infty,\infty)$, as $n\rightarrow\infty$, for independent and
weakly dependent data $\{t_{i}\}$, derived in Anevski and H\"ossjer
\cite{anevski:hossjer:2006}.
\par
The kernel estimator can be written $x_n=x_{b,n}+v_n$ with
\begin{eqnarray}
        x_{n}(t)&=&h^{-1}\int k'(u)F_{n}(t-hu)\,du,\nonumber\\
        x_{b,n}(t)&=&h^{-1}\int k'(u)F(t-hu)\,du,\label{def:xbn-dens}\\
        v_{n}(t)&=&h^{-1}\int k'(u)F_{n}^{0}(t-hu)\,du.\nonumber
\end{eqnarray}
Rescaling is done on a scale $d_n$ that is of the same asymptotic
order as $h$, so that we put $d_n=h$. The rescaled process is
\begin{eqnarray*}
   \tilde{v}_{n}(s;t_0)&=&c_{n}\int k'(u)w_{n,d_{n}}(s-u;t_{0})\,du,
\end{eqnarray*}
with $c_{n}=d_{n}^{-1}(nh)^{-1}\sigma_{n,d_{n}}$.
\par
%  ---- Limit theorem for the density  --------
\bth\label{thm:denslim} Let $\{t_{i}\}_{i\geq1}$ be a stationary
sequence with a monotone mar\-gi\-nal density function $f$ such that
$\sup_{t\in I_{t_0}}f'(t)<0$ and $f\in C^{1}(I_{t_0})$ for an open
interval $I_{t_0}\ni t_0$ where $t_{0}>0$. Assume that ${\mathbb
E}t_i^5 < \infty$. Let $x_{n}$ be the kernel density function
defined above, with $k$ a bounded and compactly supported density
such that $k'$ is bounded. Suppose that one of the following
conditions holds:
\begin{eqnarray*}
{\bf [a]}&&\{t_{i}\}_{i\geq1}\mbox{ is an i.i.d. sequence,}\\
  {\bf [b]}&&1)\;\{t_{i}\}_{i\geq1}\mbox{ is a stationary  $\phi$-mixing  sequence
 with $\sum_{i=1}^{\infty}
\phi^{1/2}(i)<\infty$} \; ; \\
&&2)\;  \mbox{$f(t_{0})=F'(t_{0})$ exists, as well as the joint
density
$f_{k}(s_{1},s_{2})$ of}\\
&&\mbox{$(t_{1},t_{1+k})$ on
$[t_{0}-\delta,t_{0}+\delta]^{2}$ for some $\delta>0$, and $k\geq
1$} \; ; \\
&&3)\; \sum_{k=1}^{\infty}M_{k}<\infty \mbox{ holds, for
}M_{k}=\sup_{t_{0}-\delta\leq s_{1},s_{2}\leq t_{0}+\delta}
|f_{k}(s_{1},s_{2})-f(s_{1})f(s_{2})|.
    \end{eqnarray*}
Then choosing $h=an^{-1/3}$ and $a>0$ an arbitrary constant, we
obtain
\begin{eqnarray*}
    n^{1/3}\{\hat{f}_n(t_{0})-f(t_{0})\}&\stackrel{\cal L}{\rightarrow}&
    aT[f'(t_{0})\cdot +\tilde{v}(\cdot ;t_{0})](0)+f'(t_0)\,a\int u k(u)\,du,
\end{eqnarray*}
as $n\rightarrow\infty$, where $\tilde{v}(s;t)$ is as in
$(\ref{eq:PaAtildev5.2.1})$, with $c=a^{-3/2}f(t_{0})^{1/2}$, and
$w$ a standard two sided Brownian motion.
 \eth
 \prf If $k'$ is bounded and
$k$ has compact support, the continuity of the map
\begin{eqnarray*}
 C(-\infty,\infty)  \ni z(s)&\mapsto&\int z(s-u)k'(u)\,du \in C(-\infty,\infty)
\end{eqnarray*}
implies that, choosing $d_{n}$ such that $c_{n}\rightarrow c$ as
$n\rightarrow\infty$ for some constant $c$, one gets:
\begin{eqnarray}
        \tilde{v}_{n}(s;t_{0})&\stackrel{\cal L}{\rightarrow}&
        c\int k'(u) w(s-u ; t_0)\,du=:\tilde{v}(s;t_{0}), \label{eq:PaAtildev5.2.1}
\end{eqnarray}
on $C(-\infty,\infty)$, as $n\rightarrow\infty$, thanks to the
continuous mapping theorem. Here $w$ is the weak limit of
$\{w_{n}\}$.  Theorems 7 and 8 of Anevski and H\"ossjer
\cite{anevski:hossjer:2006} state that
 $w_{n,\delta_n}(s,t_0)\stackrel{\cal L}{\rightarrow} B(s)$ on
$D(-\infty,\infty)$ under the respective assumptions in ${\bf [a]}$
and ${\bf [b]}$, where $B(s)$ is a two sided standard Brownian
motion. This establishes the first part of Assumption 1 for both
cases ${\bf [a]}$ and ${\bf [b]}$. Next notice that
$x_{b,n}(t)=h^{-1}\int k(\frac{t-u}{h})f(u)\,du$ is monotone. A
change of variable and a Taylor expansion in $x_{b,n}$ prove  the
second part of Assumption 1 with $A=f'(t_0)$ and
\begin{eqnarray*}
   d_n^{-1}\{x_{b,n}(t_0)-f(t_0)\}&\rightarrow& f'(t_0) \int uk(u)\,du =\Delta.
\end{eqnarray*}
\par
The statement of Assumption 2 is relegated to the appendix, see
Corollary \ref{cor:dens} in Appendix~\ref{append-dens}.
Theorem~\ref{thm:denslim} therefore holds as an application of
Theorem~\ref{thm:asympt-distr}.
\par
Let us finally check that the scale $d_n$ can be chosen so that $c_n
\to c$, as assumed at the beginning of the proof:
\begin{itemize}
\item[-] Independent data case ${\bf [a]}$: We have $\sigma_{n,d_{n}}^{2}\sim
nd_{n}f(t_{0})$, so that
\begin{eqnarray*}
        d_{n}^{-1}(nh)^{-1}\sigma_{n,d_{n}}&\sim& d_{n}^{-3/2}n^{-1/2}f(t_{0})^{1/2}.
\end{eqnarray*}
Choosing $d_{n}=a n^{-1/3}$ we get $c=a^{-3/2}f(t_{0})^{1/2}$.
\item[-] Mixing data case ${\bf [b]}$: Similar to the proof of case ${\bf [a]}$.
\end{itemize}
\vskip -6mm \eop
\bnote The present estimator was first proposed for independent data
by Fou\-g\`e\-res \cite{fougeres:1997}, who stated the strong
consistency uniformly over ${\mathbb R^+}$ for $T(f_n)$  and derived
some partial results for the limit distribution. The results for the
monotone density function estimator are similar to the results for
the Grenander estimator (the NPMLE) of a monotone density, in that
we have cube root asymptotics and a limit random variable that is a
nonlinear functional of a Gaussian process, for independent and weak
dependent data; see Prakasa Rao \cite{prakasarao:1969} and Wright
\cite{wright:1981} for the independent data cases, and Anevski and
H\"ossjer \cite{anevski:hossjer:2006} for the weak dependent data
cases. In our case however we obtain one extra term that arises from
the bias in the kernel estimator. Our estimator is really closer in
spirit to the estimator obtained by projecting the kernel estimator
on the space of monotone functions (i.e. kernel estimation followed
by isotonic regression) first proposed by Anevski and H\"ossjer
\cite{anevski:hossjer:2006}; note that we obtain the same bias term
as in Anevski and H\"ossjer \cite{anevski:hossjer:2006}. \enote
\bnote
 The results for the  long range dependence case is similar to the result
  for the isotonic regression of a kernel estimator, cf. Anevski and
  H\"ossjer \cite{anevski:hossjer:2006}.  In this situation  $\tilde{v}_{n}(s;t_0)$ is
  asymptotically a linear function  of $s$ with a random slope, implying that
  the monotone rearrangement of $g_n+\tilde{v}_n$ is just $g_n+\tilde{v}_n$
  which evaluated at zero is zero. This is due to the fact that for long range
  dependent data the limit process of the empirical process is a deterministic
  function multiplied by a random variable, cf. the remark after Theorem 12 in Anevski and H\"ossjer
  \cite{anevski:hossjer:2006}. Thus the limit distribution for the final
  estimator for long range dependent data is the same as the limit distribution
  for the kernel estimator itself, i.e. $n^{d/2}\{\hat{f}_n(t)-f(t)\}$ and
  $n^{d/2}\{f_n(t)-f(t)\}$ have the same distributional limit. See
  Cs\"org\"o and Mielniczuk \cite{csorgo:mielniczuk:1995} for a derivation of this
   limit distribution. \enote

\section{Conclusions}
We considered the feature of estimating an arbitrary monotone
function $x$, via a monotone rearrangement of a "preliminary"
estimator $x_n$ of the unknown $x$. We derived consistency and limit
distribution results for the monotonized estimator that hold under
rather general dependence assumptions.
\par
Our approach is similar in spirit to the general methods studied in
Anevski and H\"ossjer \cite {anevski:hossjer:2006} and first
introduced in the regression estimation setting  by Mammen
\cite{mammen:1991}: Start with a preliminary estimator
%(in Mammen's
%paper explicitly stated to be a kernel estimator, whereas Anevski
%and H\"ossjer's paper uses an approach that is similar to ours)
and make it monotone by projecting it on the space of monotone
functions.
%\par
The present approach can however at some point be considered
preferable: The monotone rearrangement, being basically a sorting,
is a simpler procedure than an ${\mathbb L}^2$-projection.
Furthermore the consistency and limit distribution results indicate
similar properties to Mammen's and Anevski and H\"ossjer's
estimators.  Besides, an important advantage of our estimator is the
finite sample behavior: Mammen's estimator is monotone but not
necessarily smooth; Mammen actually  studied two approaches, one
with kernel smoothing followed by monotonization and the other
approach the other way around, i.e. monotonization followed by
kernel smoothing. Mammen showed that the two proposals are first-order
equivalent. However, their finite sample size properties are very
different: the first resulting estimator is monotone but not
necessarily smooth, while the other is smooth but not necessarily
monotone, so that one needs to choose which property is more
important. This is not the case with our estimator, since if we
start with a smooth estimator of the function, e.g. a kernel
estimator, the monotone rearrangement will be smooth as well. This
can however become a disadvantage for instance when the estimand is
discontinuous: then the monotone rearrangement will "oversmooth"
since it will give a continuous result, while Mammen's estimator
will keep more of the discontinuity intact.

Some simulation studies are available in the literature, which
exhibit the small sample size behavior of the rearrangement of a
kernel estimator of a density, and compare it to different
competitors. See e.g. Foug\`eres \cite{fougeres:1997}, Meyer and
Woodroofe \cite{meyer:woodroofe:2004}, Hall and Kang
\cite{hall:kang:2005}, Chernozhukov et al.
\cite{chernozhukov:fernandezval:galichon:2007a}. These references
deal with independent data. A larger panel of dependence situations
in the comparisons would clearly be of interest, and this will be
the object of future work.

\par
Note that our results are geared towards local estimates, i.e.
estimates that use only a subset of the data and that are usually
estimators of estimands that can be expressed as non-differentiable
maps of the distribution function such as e.g. density functions,
regression functions, or spectral density functions. This differs
from global estimates, as those considered for example by
Chernozhukov et al.\cite{chernozhukov:fernandezval:galichon:2007b}
for quantile estimation.
\par
An approach similar to ours for local estimates is given in  Dette
et al. \cite{dette:neumeyer:pilz:2006},  using a modified version of
the Hardy-Littlewood-P\'olya monotone rearrangement: The first step
consists of calculating the upper level set function and is
identical to ours. However in the second step they use a smoothed
version of the (generalized) inverse, which avoids nonregularity
problems for the inverse map. The resulting estimator is therefore
not rate-optimal, and the limit distributions are standard Gaussian
due to the oversmoothing.
\par
Work has been done here using kernel based methods for the
preliminary estimator $x_n$ of $x$. Other methods, such as wavelet
based ones, are possible, and let emphasize that the only
assumptions  required are given in Assumptions 1 and 2.
\par
We have studied applications to density and regression function
estimation. Other estimation problems that are possible to treat
with our methods are e.g. spectral density estimation, considered by
Anevski and Soulier \cite{anevski:soulier:2007}, and deconvolution,
previously studied by van Es et al.
\cite{vanes:jongbloed:vanzuijlen:1998} and Anevski
\cite{anevski:1999}.

  \bibliography{bibMONREA}

\newpage

\appendix
\section{Maximal bounds for rescaled partial sum and empirical processes}
 In
this section we derive conditions under which Assumption 2 holds,
for the density and regression function estimation cases. Recall
that
\begin{eqnarray}
\tilde{g}_n(s)&=&d_n^{-1}\{x_{b,n}(t_0+sd_n)-x_{b,n}(t_0)\},\label{eq:gn-init}\\
\tilde{v}_n(s)&=&d_n^{-1}v_{n}(t_0+sd_n)\nonumber.
\end{eqnarray}
Since under Assumption 1
\begin{eqnarray*}
      y_n(s)-\{\tilde{g}_n(s)+\tilde{v}_n(s)\}&=& d_n^{-1}\{x_{b,n}(t_0)-x(t_0)\}\\
      &\rightarrow&\Delta,
\end{eqnarray*}
as $n\rightarrow \infty$, and $|\Delta|<\infty$, establishing
Assumption 2 for the process $\tilde{g}_n+\tilde{v}_n$ implies that
it holds also for the process $y_n=g_n+\tilde{v}_n$. Therefore it is
enough to establish Assumption 2 for  ${y}_n$ replaced by
$\tilde{g}_n+\tilde{v}_n$.
\par
Recall that for the cases that we cover the rescaled processes are of the form
\begin{eqnarray*}
     \tilde{v}_{n}(s;t_0)&=&c_n\int k'(u)z_n(s-u;t_0)\,du,
\end{eqnarray*}
with $z_n=w_{n,d_n}$ the local rescaled empirical process in the
density estimation case and $z_n=w_n$ the partial sum process in the
regression case. This implies that for the density estimation case
the support of $\tilde{v}_n$ is stochastic, since it depends on
$\max_{1\leq i\leq n} t_i$, while for the regression estimation case
it does not depend on the data $\{\epsilon_i\}$ and is as a matter
of fact compact and deterministic.
 \blem\label{lem:appendix}
Let $\supp(k)\subset [-1,1]$. Suppose that Assumption 1 holds.
Assume that $t_0$ has a neighbourhood
$I=[t_0-\epsilon,t_0+\epsilon]$ such that $\tau:=\sup_{t\in I  }
x'(t)<0$.   Suppose also that
\begin{eqnarray}
x_{b,n}'(t+sd_n)&\to& x'(t), \label{eq:xbnlimit}
\end{eqnarray}
as $n\to\infty$, for all $t\in I$.
\par
 Then $(\ref{eq:downdip1})$ and $(\ref{eq:downdip2})$ written
for $z_n=\tilde{g}_n+\tilde{v}_n$ are implied by the two results: \\
(A). For every
$\delta>0$ and $0<c<\infty$ there is a
  finite $M>0$
\begin{eqnarray*}
\liminf_{n\to \infty} P\left[\cap_{s\in (c,d_n^{-1}\epsilon)}\{
\tilde{v}_n(s)<\frac{M}{2}-\tau(s-c)\}\right]&>&1-\delta\label{eq:boundAAA2}.
\end{eqnarray*}
 (B). For every
$\delta>0$ and finite $M>0$ there is a finite $C$  such that for each $c>C$
\begin{eqnarray}
    \limsup_{n\rightarrow \infty}P\left\{\sup_{s\in d_n^{-1}(0,\ell(n))} \tilde{v}_n(s)>\frac{M}{2}
    -\tau(d_n^{-1}\epsilon-c)\right\}&<&\delta, \label{eq:CS-A2}
\end{eqnarray}
where $\ell(n)$ is  a deterministic function which satisfies either
of
\begin{eqnarray*}
     (i)&&\liminf_{n\to \infty}P\{\max_{1\leq i\leq n} t_i<\ell(n)\}= 1,~~~~~~~~~~~~
\end{eqnarray*}
or
\begin{eqnarray*}
    (ii)&&\ell(n)\equiv  \max \supp(x_n) {\mbox{~if~}} \limsup_{n\to \infty} \max \supp(x_n) \leq K<\infty.
\end{eqnarray*}

 \elem
Condition $(A)$ can be seen as boundedness on small sets (i.e. on the sets $(c,d_n^{-1}\epsilon
)$), while the conditions in $(B)$ are bounds outside of small sets; the small sets are really compact (of the form $(0,\epsilon)$) on the $t$-scale, and are increasing due to the rescaling done for the $s$-scale.
 \par
Condition $(B)(ii)$ is appropriate for the regression function estimation
 case, since then $\limsup_{n\to \infty} \max(\supp(x_n))$ is bounded by
 $1+\max(\supp(k))=2$, while for the density estimation case we will have
 to invoke the more subtle assumptions in $(B)(i)$.
 \par
\noindent \prf
 In order to show (\ref{eq:downdip1}), we first prove
that for each $\delta>0$ there is a $0<M<\infty$
 and a $0<c<\infty$ such that
 \begin{eqnarray}
       \liminf_{n\to \infty}P\{\sup_{s\geq c}(\tilde{g}_n+\tilde{v}_n)(s)<-M\}&\geq&1-\delta. \label{eq:firstmaxbound}
 \end{eqnarray}
 Let $\tilde{g}_n$ be defined in (\ref{eq:gn-init}). Consider the function
\begin{eqnarray*}
           k_{n}(s)&=&\left\{ \begin{array}{ll}
                                 \tau s,\mbox{ on }  (-\epsilon d_n^{-1},\epsilon d_n^{-1}),\\
                                \tau \epsilon d_n^{-1} \mbox{ on }(\epsilon d_n^{-1},\infty),\\
                                -\tau \epsilon d_n^{-1} \mbox{ on } (-\infty,-\epsilon d_n^{-1}).
                                         \end{array} \right .
\end{eqnarray*}
Then from (\ref{eq:xbnlimit}) we obtain
\begin{eqnarray*}
   \tilde{g}_n(s)&\leq & k_n(s)\mbox{ on }{\mathbb R}^+,\\
   \tilde{g}_n(s)&\geq &k_n(s)\mbox{ on }{\mathbb R}^-,
\end{eqnarray*}
for all large enough $n$, since $\tilde{g}_n$ is decreasing (as
weighted mean of decreasing functions) and  $\tilde{g}_n(0)=0$.
\par
Let $\delta$ be given and suppose part (A) of the assumptions is
satisfied, with some $M$ and arbitrary $0<c<\infty$. We will
consider the hypotheses $(B)(i)$ and $(B)(ii)$ separately:
\vskip 3mm $(B)(i)$ Since the kernel $k$ has support in $[-1,1]$ one
has $\supp(x_n)\subset(\min_{1\leq i\leq n}t_i-h,\max_{1\leq i\leq
n} t_i+h)$. Using the rescaling $t=t_0+sd_n$ this implies that
\begin{eqnarray*}
    \supp(\tilde{g}_n),\supp(\tilde{v}_n)&\subset& -t_0+d_n^{-1}(\min t_i -h,\max t_i +h)=:I_n^{(i)}.
\end{eqnarray*}
Since $t_0>\min t_i$ and
$h$ is positive, the supremum over all $s\in I_n^{(i)}$ can be
replaced by a supremum over all $s\in(c,d_n^{-1}\max t_i)$, as $n$
tends to $\infty$, and thus we need to show
\begin{eqnarray}
     \liminf_{n\rightarrow \infty} P\left\{\sup_{(c,d_n^{-1}\max t_i)} (\tilde{g}_n+\tilde{v}_n)
     (s)<-M\right\}&\geq&1-\delta. \label{eq:boundAAA}
\end{eqnarray}
Then for $c\geq 3M/2|\tau|$, we will have $k_n(c)= -3M/2$. This
implies that for $c\geq 3M/2|\tau|$,
\begin{eqnarray*}
      P\left\{\sup_{(c,d_n^{-1}\max t_i)} y_n(s)<-M\right\}
&\geq&   P (\cap_{s\in (c,d_n^{-1}\epsilon)}\{ \tilde{v}_n(s)<\frac{M}{2}-\tau(s-c)\}\\
 &&
\cap \{\sup_{s\in d_n^{-1}(\epsilon,\max t_i)}
\tilde{v}_n(s)<\frac{M}{2}-\tau(d_n^{-1}\epsilon-c)\}),
\end{eqnarray*}
so that $(\ref{eq:boundAAA})$ follows from the two results
\begin{eqnarray}
\liminf_{n\to \infty} P\left\{\sup_{s\in d_n^{-1}(\epsilon,\max
t_i)} \tilde{v}_n(s)<\frac{M}{2}-\tau(d_n^{-1}\epsilon-c)\right\}
&>&1-\delta,\label{eq:boundAAA1}\\
\liminf_{n\to \infty} P\left[\cap_{s\in (c,d_n^{-1}\epsilon)}\{
\tilde{v}_n(s)<\frac{M}{2}-\tau(s-c)\}\right]&>&1-\delta\label{eq:boundAAA2}.
\end{eqnarray}
\par
The relation $(\ref{eq:boundAAA2})$ is satisfied by assumption $(A)$ and thus we need to treat $(\ref{eq:boundAAA1})$.  Let $\ell$ be the deterministic function given in assumption $(B)(i)$. Note first that
\begin{eqnarray*}
   &&P\left\{\sup_{d_n^{-1}(\epsilon,\ell(n))}
\tilde{v}_n(s)<\frac{M}{2}-\tau(d_n^{-1}\epsilon-c)\right\}\\
&\leq&  P\left\{\sup_{d_n^{-1}(\epsilon,\ell(n))}
\tilde{v}_n(s)<\frac{M}{2}-\tau(d_n^{-1}\epsilon-c)| \max_{1\leq i\leq n} t_i<\ell(n)\right\}+   P\left\{\max_{1\leq i\leq n} t_i>\ell(n)\right\}\\
&<&P\left\{\sup_{d_n^{-1}(\epsilon,\ell(n))}
\tilde{v}_n(s)<\frac{M}{2}-\tau(d_n^{-1}\epsilon-c)| \max_{1\leq i\leq n} t_i<\ell(n)\right\}+\delta
\end{eqnarray*}
for all $n\geq N$ for some $N$,  since $\lim_{n\to \infty}P\left\{\max_{1\leq i\leq n} t_i>\ell(n)\right\}= 0$. Therefore, for all $n\geq N$, we have
\begin{eqnarray*}
&&P\left\{\sup_{d_n^{-1}(\epsilon,\max t_i)} \tilde{v}_n(s)<\frac{M}{2}-\tau(d_n^{-1}\epsilon-c)\right\}\\
&\geq &P\left\{\sup_{d_n^{-1}(\epsilon,\ell(n))}
\tilde{v}_n(s)<\frac{M}{2}-\tau(d_n^{-1}\epsilon-c)\,|\,
\max_{1\leq i\leq n} t_i<\ell(n)\right\}P\left\{\max_{1\leq i\leq n} t_i<\ell(n)\right\} \\
&\geq&\left(P\left\{\sup_{d_n^{-1}(\epsilon,\ell(n))}
\tilde{v}_n(s)<\frac{M}{2}-\tau(d_n^{-1}\epsilon-c)\right\}-\delta\right)
P\left\{\max_{1\leq i\leq n} t_i<\ell(n)\right\}\\
  &\geq&\left(P\left\{\sup_{d_n^{-1}(0,\ell(n))} \tilde{v}_n(s)<\frac{M}{2}-\tau(d_n^{-1}\epsilon-c)\right\} -\delta\right)
  P\left\{\max_{1\leq i\leq n} t_i<\ell(n)\right\}.
\end{eqnarray*}
Thus since
$\lim_{n\to\infty}P\left\{\max_{1\leq i\leq n}
t_i<\ell(n)\right\}=1$, taking complements leads to
(\ref{eq:boundAAA1}) as soon as for $c>C$
\begin{eqnarray*}
    \limsup_{n\to \infty}P\left\{\sup_{d_n^{-1}(0,\ell(n))}
    \tilde{v}_n(s)>\frac{M}{2}-\tau(d_n^{-1}\epsilon-c)\right\}&<&\delta,
%    \label{eq:CS-AAA1}
\end{eqnarray*}
i.e. (\ref{eq:CS-A2}).

\vskip 3mm $(B)(ii)$. It follows from the definition of $K$ and from
$\supp(k)\subset [-1,1]$  that $\supp (x_n)\subset (-h,K+h)$, so
that
\begin{eqnarray*}
    \supp (\tilde{g}_n),\supp(\tilde{v}_n)&\subset& -t_0+d_n^{-1}(-h,h+K)=:I_n^{(ii)}.
\end{eqnarray*}
Again this implies that the supremum of $\tilde{v}_n$ over
$I_n^{(ii)}$ can be replaced by a supremum over all $s\in
(c,d_n^{-1}K)$ and thus (\ref{eq:firstmaxbound}) will follow as soon
as
\begin{eqnarray}
   \liminf_{n\to \infty} P\left\{\sup_{(c,d_n^{-1}K)} (\tilde{g}_n+\tilde{v}_n)(s)<-M\right\}
   &>&1-\delta.\label{eq:boundAAAprim}
\end{eqnarray}
For arbitrary $M$ and $c\geq 3M/2$ we have
\begin{eqnarray*}
      && P\left\{\sup_{(c,d_n^{-1}K)}
      (\tilde{g}_n+\tilde{v}_n)(s)<-M\right\}
\\&\geq&   P (\cap_{s\in (c,d_n^{-1}\epsilon)}\{ \tilde{v}_n(s)<\frac{M}{2}-\tau(s-c)\}
\\ &&
\cap \{\sup_{s\in d_n^{-1}(\epsilon,K)}
\tilde{v}_n(s)<\frac{M}{2}-\tau(d_n^{-1}\epsilon-c)\}),
\end{eqnarray*}
so that $(\ref{eq:boundAAAprim})$ follows from
\begin{eqnarray}
\liminf_{n\to \infty} P\left\{\sup_{s\in d_n^{-1}(\epsilon,K)}
\tilde{v}_n(s)<\frac{M}{2}-\tau(d_n^{-1}\epsilon-c)\right\}&>&1-\delta,\label{eq:boundAAA1prim}
\end{eqnarray}
and $(\ref{eq:boundAAA2})$, which ends the derivation for the case $(ii)$.
\vskip 3mm
\par
Now we prove that with $M$ as above
\begin{eqnarray}
   \liminf_{n\to\infty} P\left\{\inf_{\inf I_1 \leq s\leq \sup I_0} y_n(s)\geq
   -M\right\}
   &>&1-\delta. \label{eq:secondmaxbound}
\end{eqnarray}
Note that with $M=M_c$ corresponding to the bound for $c$, we have
$k_n(c)=-3M_c/2$
and thus ${\mathbb E}y_n(c)\leq -3M_c/2\leq -M_c-M_c/2$. Since
$\tilde{g}_n(s)\to As$ on compact intervals, if $n$ is large enough
then  we have ${\mathbb E}y_n(s)\geq As-\epsilon$ for each
$\epsilon>0$ arbitrarily small. Thus for $s_M=-M_c/2A$
\begin{eqnarray*}
  {\mathbb E}y_n(s_M)&\geq&-\frac{M_c}{2}-\epsilon,
\end{eqnarray*}
for $n$ large enough. Finally from $(\ref{eq:firstmaxbound})$  and
by the symmetry of the distribution of $\tilde{v}_n$ around 0, we
have that with $M$ replaced by $\max\{M_c,-2A\sup I_0\}$, both
$(\ref{eq:firstmaxbound})$ and $(\ref{eq:secondmaxbound})$ hold, and
(\ref{eq:boundAAA1}) is proven.
\par
 Equation (\ref{eq:boundAAA2}) can be
proven in a similar way, which yields the lemma.  \eop
Lemma \ref{lem:appendix} states conditions $(A)$ and $(B)$ as
sufficient conditions for Assumption 2.
To further simplify condition $(B)$ in Lemma \ref{lem:appendix}, using Boole's
inequality and the stationarity of the process $\tilde{v}_n$ we get
in both cases $(i)$ and $(ii)$
\begin{eqnarray}
    && P\left\{\sup_{d_n^{-1}(0,\ell(n))} \tilde{v}_n(s)>\frac{M}{2}-\tau(d_n^{-1}\epsilon-c)\right\}\label{eq:probtomaj}\\
     &&\leq d_n^{-1}\ell(n) P\left\{\sup_{(0,1)}
     \tilde{v}_n(s)>\frac{M}{2}-\tau(d_n^{-1}\epsilon-c)\right\},\nonumber
%     (ii)&&P(\sup_{d_n^{-1}(0,K)} \tilde{v}_n(s)>\frac{M}{2}-\tau(d_n^{-1}\epsilon-c))\\
%     &&\leq  d_n^{-1}K P(\sup_{(0,1)} \tilde{v}_n(s)>\frac{M}{2}-\tau(d_n^{-1}\epsilon-c)).
\end{eqnarray}
where $\ell(n)$ is defined for hypothesis $(i)$ and replaced by $K$
when dealing with hypothesis $(ii)$. As a consequence, in Case $(i)$
(resp. Case $(ii)$) the probability (\ref{eq:probtomaj}) will
converge to 0 as soon as
\begin{eqnarray*}
   \alpha(n)&:=&  P\left\{\sup_{(0,1)} \tilde{v}_n(s)>\frac{M}{2}-\tau(d_n^{-1}\epsilon-c)\right\}\to 0
\end{eqnarray*}
faster than $d_n^{-1}\ell(n)\to \infty$, i.e. that
$\alpha(n)=o(d_n\, \ell(n)^{-1})$ as $n\rightarrow \infty$ (resp.
$\alpha(n)=o(d_n)$). The following conditions are thus respectively
sufficient to insure that (\ref{eq:probtomaj}) tends to 0, as
$n\rightarrow \infty$:
\begin{eqnarray*}
      (i)&&P\left\{\max_{1\leq i\leq n} t_i<\ell(n)\right\}\rightarrow
      1 {\mbox{~and~}} d_n^{-1}\ell(n)\alpha(n)\rightarrow 0,\\
      (ii)&&d_n^{-1} \alpha(n)\to 0.
\end{eqnarray*}
\par
Finally, the examination of the convergence of $\alpha(n)$ can be made
in two steps via the standard partition
 \begin{eqnarray}
\alpha(n) &\leq& P\left[\sup_{s,s'\in(0,1)}| \tilde{v}_n(s)-\tilde{v}_n(s')|>\frac{1}{2}
 \left\{\frac{M}{2}-\tau(d_n^{-1}\epsilon-c)\right\}\right]\nonumber \\
 &&+ \; P\left[\tilde{v}_n(0)>\frac{1}{2}\left\{\frac{M}{2}-\tau(d_n^{-1}\epsilon-c)\right\}\right].\label{eq:maxpartition}
 \end{eqnarray}
In the sequel we will bound the two terms of the right-hand side of $(\ref{eq:maxpartition})$ separately, for the
  density and regression estimation problems treated in this paper: See subsections A.1 and A.2.
\vskip 3mm To further simplify $(A)$ in Lemma \ref{lem:appendix},
note that
\begin{eqnarray*}
  &&\cap_{s\in (c,d_n^{-1}\epsilon)} \{   \tilde{v}_n(s)<\frac{M}{2}-\tau(s-c)\}\\
    &\supset & \cap_{i\in {\mathbb Z}\cap (c,d_n^{-1}\epsilon)} \{ \sup_{s\in [i,i+1)}
     \tilde{v}_n(s)<\frac{M}{2}-\tau(i-c) \}=:A_n.
\end{eqnarray*}
Thus, taking complements, part $(A)$ of Lemma \ref{lem:appendix}
follows as soon as for every $\delta$ and arbitrary $0<c<\infty$
there is a $0<M<\infty$ such that $\limsup_{n\to\infty}
P(A_n^c)<\delta$. However,
\begin{eqnarray*}
    P(A_n^c)&\leq& \sum_{i\in {\mathbb Z}\cap (c,d_n^{-1}\epsilon)}
    P\left\{\sup_{s\in[i,i+1)} \tilde{v}_n(s)> \frac{M}{2}-\tau(i-c)\right\}\\
&{\leq}& \sum_{i=[c]}^{[d_n^{-1}\epsilon]} P\left\{\sup_{s\in [0,1)}
\tilde{v}_n(s)> \frac{M}{2}-\tau(i-c)\right\},
\end{eqnarray*}
where the equality follows from the stationarity of $\tilde{v}_n$.
In the sequel we will establish maximal inequalities of the form
\begin{eqnarray}
 P\left\{\sup_{s\in [0,1)} \tilde{v}_n(s)> a\right\}&\leq &C´ a^{-p} \label{eq:maxinequal}
\end{eqnarray}
 for some constant $p>1$; assume for now that these are established. Then
\begin{eqnarray*}
   \limsup_{n\to \infty}P(A_n^c)&\leq & \sum_{i=[c]}^{\infty} C´ \frac{1}{(\frac{M}{2}-\tau(i-c)))^p} \\
&\leq & \frac{C}{p|\tau|^p}\left(\frac{2}{M}\right)^{p-1}\\
&<&\delta,
\end{eqnarray*}
where the next to last inequality holds by an integral approximation
of the series and the last by choosing $M=M(\delta)>2(C
/p\delta|\tau|^p)^{1/(p-1)}$. Thus assumption $(A)$ in Lemma
\ref{lem:appendix} follows from $(\ref{eq:maxinequal})$ with $p>1$;
inequalities of the form $(\ref{eq:maxinequal})$ will next be
treated.

\subsection{Maximal bounds for the rescaled partial sum process}\label{append-reg}
Let $k$ be a kernel which is bounded, piecewise differentiable, with
a bounded derivative, say $0 \leq |k'| \leq \alpha$. Assume that the
sequence $h=h_n$ is such that $nh \to \infty$. We have (see
(\ref{eq:vnrescaled}))
\begin{eqnarray*}
        \tilde{v}_n(s,t_{0})&\stackrel{\cal L}{=}& d_n^{-1}(nh)^{-1} \sigma_{\hat{n}}
        \int \bar{w}_{\hat{n}}(s-u)k'(u)\,du,
\end{eqnarray*}
where $d_n=h$ is chosen so that $d_n^{-1}(nh)^{-1} \sigma_{\hat{n}}
\rightarrow 1$ and $\hat{n}=[nh]$. Now $\bar{w}_{\hat{n}}$ is
asymptotically equivalent to the piecewise constant partial sum
process which we therefore will use for notational simplicity, and
which we denote (with a slight abuse of notation) with
$w_{\hat{n}}$.
\par
We show the convergence of $\alpha(n)$ in which $\tilde{v}_n(s)$ is
replaced by  $w_{\hat{n}}(s)$: this will be sufficient since
\begin{eqnarray*}
|\tilde{v}_{n}(s)|&\leq& \sup_{u \in [-1,1]}|w_{\hat{n}}(s-u)|  \int
|k'(u)|\,du,
\end{eqnarray*}
and thus
\begin{eqnarray*}
   \sup_{s\in (0,1)}  |\tilde{v}_{n}(s)|&\leq&     c\,\sup_{s\in (0,1)} \sup_{u\in [-1,1]}|w_{\hat{n}}(s-u)|   \\
   &\leq &c\,\sup_{s\in[-1,2]} |w_{\hat{n}}(s)|,
\end{eqnarray*}
   with $c=\int |k'(u)|\,du$, and since the behaviour of the process $w_{\hat{n}}$
   on $(0,1)$ and on $(-1,2)$ is qualitatively the same.

  \bpr \label{prop:maj}
  Let $p\geq 2$ be given and assume that the sequence $\{\epsilon_i\}_{i\geq 1}$
  satisfies $\max({\mathbb E}\epsilon_1^2,{\mathbb E}\epsilon_1^p)<\infty$. Then under the assumptions of Theorem 5
\begin{eqnarray*}
P\left( \sup_{s \in (0,1)} {w}_{\hat{n}}(s) > M/2 - \tau(\delta
h^{-1} -C) \right) &\leq& Ch^p,
\end{eqnarray*}
where $C$ is a finite constant.
 \epr
\prf %
Let $a:=M/2 + \tau c $ and $b:= - \tau \delta$. In a first step, we
obtain a majoration of
\begin{eqnarray*}
P\left( \sup_{s,s' \in (0,1)} | {w}_{\hat{n}}(s) -
{w}_{\hat{n}}(s') |
>a +b h^{-1} \right)
\end{eqnarray*}
in the 3 dependence situations listed in Theorem \ref{thm:reglim}.
One has
\begin{eqnarray}
     w_{\hat{n}}(s)-w_{\hat{n}}(s') &=&\frac{1}{\sigma_{\hat{n}}} S_n(s,s'), \label{eq:wn_diff}
\end{eqnarray}
with
\begin{eqnarray*}
        S_n(s,s')&=&\sum_{i=[s'\hat{n}]+1}^{[s\hat{n}]} \epsilon_i.
\end{eqnarray*}
\vskip 3mm {\bf [a]} If $\{\epsilon_i\}$ is an i.i.d. sequence the
moment bound in Theorem 2.9 in Petrov \cite{petrov:1994} implies
\begin{eqnarray*}
       {\mathbb E}|S_n(s,s')|^p&\leq & c(p)\left(    \sum_{i=[s'\hat{n}]+1}^{[s\hat{n}]}
        {\mathbb E} |\epsilon_i|^p  +\left(\sum_{i=[s'\hat{n}]+1}^{[s\hat{n}]}
        {\mathbb E} (\epsilon_i)^2\right)^{p/2}\right) \\
       &\leq &c' \left(   |s-s'| \hat{n}+|s-s'|^{p/2} \hat{n}^{p/2}\right)
       \end{eqnarray*}
where $c(p)$ depends on  $p$ only and $c'=c(p)\cdot
\max(||\epsilon_1||_2,{\mathbb E}|\epsilon_1|^p)$.
\vskip 3mm {\bf [b]} If $\{\epsilon_i\}$ is a stationary sequence
that is $\alpha$-mixing (and thus also $\phi$-mixing) satisfying the
mixing condition $(\ref{eq:mixing})$, then Theorem 1 in Doukhan
\cite{doukhan:1994} implies
\begin{eqnarray*}
       {\mathbb E}|S_n(s,s')|^p&\leq &\max(\hat{n}|s-s'|M_{p,\epsilon},\hat{n}^{p/2}|s-s'|^{p/2}M_{p,2}^{p/2}),
\end{eqnarray*}
where $M_{p,\epsilon}=||\epsilon_i||^{p}_{p+\epsilon}$, and thus
\begin{eqnarray*}
      {\mathbb E}|S_n(s,s')|^p&\leq &c''\max(\hat{n} |s-s'|^p, \hat{n}^{p/2}|s-s'|^2 ),
\end{eqnarray*}
with $c''=\max(({\mathbb E} |\epsilon_i |^{p+\epsilon}
)^{p/(p+\epsilon)},({\mathbb E} |\epsilon_i |^{2+\epsilon}
)^{2/(2+\epsilon)})$.
\par
Therefore, for both independence and weak dependence cases, equation
(12.42) of Billingsley \cite{billingsley:1968} is satisfied, so that
 Theorem 12.2 in Billingsley implies
\begin{eqnarray*}
     &&P\left(\sup_{s,s' \in (0,1)} | {w}_{\hat{n}}(s) -
{w}_{\hat{n}}(s') |
>a +b h^{-1} \right)\\
&=&P\left(\max_{k \in ([s'\hat{n}]+1,[s\hat{n}])} \sigma_{\hat{n}}^{-1} \sum_{i=[s'\hat{n}]+1}^{k} \epsilon_i>  a +b h^{-1} \right)\\
&\leq&\frac{K'_{p}C(\hat{n})}{\sigma_{\hat{n}}^p (a +b h^{-1} )^p}
\end{eqnarray*}
where $C(\hat{n})=c'(\hat{n}+\hat{n}^{p/2})$ for i.i.d. data and $C(\hat{n})=c''\max(\hat{n},\hat{n}^{p/2})$ in the mixing case. Since in both cases  $\sigma_{\hat{n}}=\hat{n}^{1/2}$ and thus $\hat{n}/\sigma_{\hat{n}}^p=\hat{n}^{1-p/2}$ and $\hat{n}^{p/2}/\sigma_{\hat{n}}^p=1$, we get the bound
\begin{eqnarray}
P\left(\sup_{s,s' \in (0,1)} | {w}_{\hat{n}}(s) -
{w}_{\hat{n}}(s') |
>a +b h^{-1} \right)&\leq& C h^{p},\label{eq:incrementmaxbound}
\end{eqnarray}
if $p\geq 2$.
\vskip 3mm {\bf [c]} In the long range dependent case we have
\begin{eqnarray*}
E(S_{\tilde{n}}^{2})&\sim&\eta_{r}^{2}l_{1}(\tilde{n})\tilde{n}^{2\beta},
\end{eqnarray*}
with $l_{1}$ as in $(\ref{eq:PaAl1})$, and according to de Haan
\cite{dehaan:1970}, equation (12.42) in Billingsley
\cite{billingsley:1968} is satisfied, with
\begin{eqnarray*}
\gamma&=&2,\\
\alpha&=&2\beta,\\
u_{l}&=&\{C_{1}\eta_{r}^{2}l_{1}(\tilde{n})\}^{1/2\beta},
\end{eqnarray*}
 for some constant $C_{1}>0$. Theorem 12.2 in Billingsley \cite{billingsley:1968} then leads to
\begin{eqnarray*}{\label{eq:PaAA}}
      P\left(\max_{k \in ([s'\hat{n}]+1,[s\hat{n}])} \sigma_{\hat{n}}^{-1} \sum_{i=[s'\hat{n}]+1}^{k}
       \epsilon_i>  a +b h^{-1} \right)&\leq
       &\frac{K'_{2,2\beta}}{(a+bh^{-1})^2
       \sigma_{\tilde{n}}^{2}}\left(\sum_{i=1}^{\tilde{n}}u_{i}\right)^{2\beta}\\
      &=&\frac{C}{(a+bh^{-1})^{2}},
\end{eqnarray*}
with $C=C_{1}K'_{2,2\beta}$, as
$\sigma_{\hat{n}}^2=\eta_{r}^{2}l_{1}(\tilde{n})\hat{n}^{2\beta}$.
Thus in the long range dependent case $(\ref{eq:incrementmaxbound})$
holds for $p=2$.
 \par
In a second step, using $w_n(s)-w_n(s')\stackrel{\cal L}{=}
w_n(s-s')$, one can then deduce from (\ref{eq:incrementmaxbound})
that
\begin{eqnarray*}
P\left[ \sup_{s \in (0,1)}  {w}_n(s)
>a +b h^{-1} \right]  &\leq& C'' h^p,
\end{eqnarray*}
where $C''>0$, which together with $(\ref{eq:incrementmaxbound})$ and
via $(\ref{eq:maxpartition})$ ends the proof.
\eop
%
%\vskip 10mm \noindent
\bcor \label{cor:reg}Suppose the assumptions of Theorem~\ref{thm:reglim} are satisfied;
then Assumption 2 holds for $y_n=g_n+{\tilde v}_n$ and for $y$ as
defined in (\ref{eq:PaAy}) in each context {\bf [a]}, {\bf [b]} and
{\bf [c]} listed in Theorem~\ref{thm:reglim}.
 \ecor
\prf Note first that if $x_{b,n}$ is defined by (\ref{def:xbn-reg}),
if $m$ is a $C^1$-function, and $k$ is a kernel with compact
support, then $x'_{b,n}(t+s d_n) \to m'(t)$ for each $t$ when $n \to
\infty$. Besides, a consequence of Proposition \ref{prop:maj} is
that
\begin{eqnarray*}
 \limsup_{n \to \infty} P\left[  \sup_{s \in (0,h^{-1}K)}
\tilde{v}_n(s) < \frac{M}{2} - \tau(h^{-1}\epsilon -c)  \right]
&\geq &1- \delta.
\end{eqnarray*}
Thus,  condition $(B)(ii)$ in Lemma \ref{lem:appendix} is satisfied
as soon as
%\begin{eqnarray*}
         $ d_n^{-1}d_n^{p}=d_n^{p-1}\rightarrow0,$ which is equivalent to $p>1$.
%\end{eqnarray*}
Thus the existence of two moments suffices to get condition (B) (ii)
of Lemma~\ref{lem:appendix} for i.i.d., mixing and subordinated
Gaussian long range dependent sequences. Condition $(A)$ of Lemma
$\ref{lem:appendix}$ follows immediately from Proposition
\ref{prop:maj}. Hence Lemma~\ref{lem:appendix} can be applied, so
that Assumption 2 holds for $y_n=g_n+{\tilde v}_n$. An analogous result
for $y$ defined by (\ref{eq:PaAy}) follows easily from the
stationarity and finite second moments of ${\tilde v}_n$.
\eop
\subsection{Maximal bounds for the rescaled empirical process}\label{append-dens}
The rescaled process is
\begin{eqnarray*}
\tilde{v}_{n}(s;t_0)&=&c_{n}\int
k'(u)w_{n,d_{n}}(s-u;t_{0})\,du
\end{eqnarray*}
 with
$c_{n}=d_{n}^{-1}(nh)^{-1}\sigma_{n,d_{n}}$. Note that, similarly to
the regression case, deriving the maximal bound for the process
$w_{n,d_{n}}$ implies the maximal bound for the process
$\tilde{v}_n$.
\bpr \label{prop:maj-dens}  Under the assumptions of Theorem
\ref{thm:denslim}, there exists a positive constant $C'$ such that
\begin{eqnarray*}
P\left( \sup_{s \in (0,1)} {w}_{n,d_n}(s) > M/2 - \tau(\epsilon
h^{-1} -c) \right) &\leq& C' h^5,
\end{eqnarray*}
\epr
\prf
 Let $a:=M/2 + \tau c $ and $b:=
- \tau \epsilon$. For independent or mixing data satisfying the assumptions of Theorem
\ref{thm:denslim}, Lemma C2 in Anevski and H\"ossjer
\cite{anevski:hossjer:2006} implies that
\begin{eqnarray}
    &&P\left\{\sup_{s,s'\in  (0,1)}|w_{n,d_n}(s;t_0)-w_{n,d_n}(s';t_0)|\geq a+bh^{-1}
\right\}\nonumber\\
&&\leq
\frac{K}{(a+bh^{-1})^4+(a+bh^{-1})^5}\label{eq:momentbound1iid}
\end{eqnarray}
for a positive constant $K$.
 \par
 From $(\ref{eq:momentbound1iid})$  one can deduce the corresponding bound for
 \begin{eqnarray*}
 &&P\left\{\sup_{s (0,1)}|w_{n,d_n}(s;t_0)|\geq a+bh^{-1}
\right\}
\end{eqnarray*}
which together with $(\ref{eq:momentbound1iid})$  and via
$(\ref{eq:maxpartition})$, implies the statement of the proposition.
\eop
\bcor \label{cor:dens}Suppose the assumptions of Theorem~\ref{thm:denslim} are
satisfied; then Assumption 2 holds for $y_n=g_n+{\tilde v}_n$ and
for $y$ as defined in (\ref{eq:PaAy}) in each context {\bf [a]} and
{\bf [b]} listed in Theorem~\ref{thm:denslim}.
 \ecor
\prf Note first that if $x_{b,n}$ is defined by
(\ref{def:xbn-dens}), if $f$ is a $C^1$-function, and $k$ is a
kernel with compact support, then $x'_{b,n}(t+s d_n) \to f'(t)$ for
each $t$ when $n \to \infty$. Besides, a consequence of Proposition
\ref{prop:maj-dens} is that
\begin{eqnarray*}
 \liminf_{n \to \infty} P\left(  \sup_{s \in (0,h^{-1}\max t_i)}
\tilde{v}_n(s) < \frac{M}{2} - \tau(h^{-1}\epsilon -c)  \right)
&\geq &1- \delta
\end{eqnarray*}
if there exists  a function $\ell(n)$  such that
\begin{eqnarray}\label{eq:last}
P\{\max_{1\leq i\leq n} t_i \leq \ell(n)\} \to 1 {\mbox{~and~}}
h^{p-1} \ell(n)\to 0,
\end{eqnarray}
as $n\to \infty$, and with $p\geq 5$.  Note that
\begin{eqnarray*}
 P\left\{ \max_{1\leq i\leq n} t_i > \ell(n)\right\}&\leq& n \frac{{\mathbb E} |
t_i |^p}{\ell(n)^p} = \frac{nI}{\ell(n)^p}.
\end{eqnarray*}
Thus the conditions in (\ref{eq:last}) are implied by
%\begin{eqnarray*}
 $$         \ell(n)^{-p}n\rightarrow
 0,\;d_n^{-1}\ell(n)d_n^{p}\rightarrow0,$$
 which is equivalent to
$$          n^{1/p}<<\ell(n)<<d_n^{1-p}.$$
%\end{eqnarray*}
In the case of i.i.d. and mixing data, when $d_n=n^{-1/3}$, $p$
should satisfy
\begin{eqnarray*}
       \frac{1}{p}<\frac{p-1}{3}&\Leftrightarrow& (p^2-p)>3\\
       &\Leftrightarrow& p>\frac{1+\sqrt{13}}{2}.
\end{eqnarray*}
This, together with the restriction $p\geq 5$ in Proposition
\ref{prop:maj-dens}, implies that the existence of five moments
suffices to establish $(B)(i)$ in Lemma \ref{lem:appendix} for
i.i.d. and mixing data. Condition $(A)$ in Lemma \ref{lem:appendix}
is immediate from Proposition \ref{prop:maj-dens}. Hence Lemma
\ref{lem:appendix} can be applied, so that Assumption 2 holds for
$y_n=g_n+{\tilde v}_n$. An analogous result for $y$ defined by
(\ref{eq:PaAy}) follows easily from the stationarity and finite
second moments of ${\tilde v}_n$.
\eop

\end{document}